\documentclass[12pt]{amsart}
\usepackage{graphicx}
\usepackage{latexsym,amscd,amssymb,epsfig}

\numberwithin{equation}{section}

\newcommand{\aaf}{\mathfrak a}

\newcommand{\N}{\mathbb{N}}

\newcommand{\BarRes}{{\bf \sf B}}
\newcommand{\NBarRes}{{\bf \sf NB}}
\newcommand{\CRes}{{\bf \sf C}}
\newcommand{\FRes}{{\bf \sf F}}
\newcommand{\HComplex}{{\bf \sf HC}}

\newcommand{\Bc}{\mathcal{B}}
\newcommand{\Gc}{\mathcal{G}}
\newcommand{\Ec}{\mathcal{E}}
\newcommand{\Rc}{\mathcal{R}}
\newcommand{\Lc}{\mathcal{L}}
\newcommand{\SM}{\rm{MinGen}}
\newcommand{\Mmc}{\mathcal{M}}

\newcommand{\wb}{{\bf w}}

\DeclareMathOperator{\Hilb}{Hilb}
\DeclareMathOperator{\id}{id}
\DeclareMathOperator{\ini}{in_{\prec}}

\DeclareMathOperator{\lex}{lex}

\DeclareMathOperator{\op}{op}

\DeclareMathOperator{\supp}{supp}
\DeclareMathOperator{\Tor}{Tor}

\newcommand{\on}[1]{\operatorname{#1}}
\newcommand{\tol}{\longrightarrow}

\newcommand{\poinr}[2]{P^{#2}_{#1}(\ul{x},t)}
\newcommand{\hilb}[1]{\Hilb_{#1}(\ul{x},t)}

\newcommand{\Gu}{\Gamma_\uparrow}
\newcommand{\Gd}{\Gamma_\downarrow}

\DeclareMathOperator{\subg}{\unlhd}

\DeclareMathOperator{\tensor}{\otimes}
\DeclareMathOperator{\iso}{\cong}

\DeclareMathOperator{\Dirsum}{\bigoplus}

\DeclareMathOperator{\pnt}{\raise 0.5mm \hbox{\large\bf.}}

\theoremstyle{plain}
\newtheorem{Thm}{\bf Theorem}[section]
\newtheorem{Lem}[Thm]{\bf Lemma}
\newtheorem{Cor}[Thm]{\bf Corollary}
\newtheorem{Prop}[Thm]{\bf Proposition}

\theoremstyle{definition}
\newtheorem{Rem}[Thm]{\bf Remark}
\newtheorem{Ex}[Thm]{\bf Example}

\newtheorem{Def}[Thm]{\bf Definition}

\newcommand{\D}{\displaystyle}

\newcommand{\ul}[1]{\underline{#1}}

\newcommand{\disbreak}{\allowdisplaybreaks}

\begin{document}

\title[Resolutions via discrete Morse theory]{Resolution of the residue
class field via algebraic discrete Morse theory}

\date{\today}

\author{Michael J\"ollenbeck}
\address{Fachbereich Mathematik und Informatik\\
Philipps-Universit\"at Marburg\\
35032 Marburg, Germany}
\email{joella@mathematik.uni-marburg.de}

\author{Volkmar Welker}
\address{Fachbereich Mathematik und Informatik\\
Philipps-Universit\"at Marburg\\
35032 Marburg, Germany}
\email{welker@mathematik.uni-marburg.de}

\thanks{Both authors were supported by EU
      Research Training Network
``Algebraic Combinatorics in Europe'', grant HPRN-CT-2001-00272.}

%------------------------------------------------------------------------------
%
% Abstract
%
%------------------------------------------------------------------------------
\begin{abstract}
Forman's Discrete Morse theory is studied from an algebraic
viewpoint. Analogous to independent work of Emil Sk\"old\-berg we show,
that this theory can be extended to chain complexes of free modules over a
ring.  We provide three applications of this theory:

\begin{itemize}
\item[(i)] We construct a new resolution of the residue class field
  $k$ over the $k$-algebra $A$, where $A=k[x_1, \ldots, x_n]/\aaf$ is the
  quotient of the commutative polynomial ring in $n$ indeterminates by an
  ideal $\aaf$. This resolution is a commutative analogue of the Anick
  resolution in the non-commutative case.
  We prove minimality of the resolution if $\aaf$ admits a quadratic Gr\"obner
  basis or if $\ini(\aaf)$ is a complete intersection.
\item[(ii)] Now let $A = k\langle x_1, \ldots, x_n \rangle / \aaf$ be the
  quotient of the polynomial ring in $n$ non-commuting indeterminates by a
  two-sided ideal $\aaf$.
  Sk\"oldberg shows how to construct the Anick resolution of $A$ as well as
  the two-side Anick resolution
  via Algebraic Discrete Morse theory. We derive the same result and prove in addition
  the minimality of these resolutions and the rationality of the
  Poincar\'e-Betti series in special cases.
\item[(iii)] In the situation of (ii) we construct a resolution of $A$ as an
  $A\otimes A^{\op}$-module. We show that this resolution is minimal in
  special cases and thereby generalize a result by BACH used to calculate
  Hochschild homology in theses cases.
\end{itemize}
\end{abstract}

\maketitle
%------------------------------------------------------------------------------
%
% Introduction
%
%------------------------------------------------------------------------------
\section{Introduction}\label{0}

Discrete Morse theory as developed by Forman
\cite{forman1},\cite{forman2} allows to construct, starting from a
(regular) CW-complex, a new homotopy equivalent CW-complex with
fewer cells.  In this paper we describe and apply an algebraic
version of this theory, which we call `Algebraic Discrete Morse
Theory.' An analogous theory was developed by Sk\"oldberg
\cite{skoed}. We consider chain complexes $C_\bullet = (C_i,
\partial_i)_{i \geq 0}$ of free modules $C_i$ over a ring $R$. A
priori we alway fix a basis the $C_i$ -- the basis elements play
the role of the cells in the topological situation. Then applying
Algebraic Discrete Morse theory constructs a new chain complex of
free $R$-modules such that the homology of the two complexes
coincides.

In Section \ref{secmorse} we describe Algebraic Discrete Morse theory
and apply the theory in the remaining sections in three situations.

In Section \ref{commcase} we consider
resolutions of the field $k$ over a quotient $A = S/\aaf$ of the
commutative polynomial ring
$S= k[x_1, \ldots, x_n]$ in $n$ variables by an ideal $\aaf$.
We construct a free resolution of $k$ as an $A$-module, which can be seen as a
generalization of the Anick resolution to the commutative case.
Our resolution is minimal, if $\aaf$ admits a quadratic Gr\"obner basis.
Also we give an explicit description of the minimal resolution of $k$,
if the initial ideal of $\aaf$ is a complete intersection.

Section \ref{noncommcase} considers the same situation in the non-commutative
case. We apply Algebraic Discrete Morse theory in order to
obtain the Anick resolution of the residue class field $k$ over
$A = k\langle x_1,\ldots, x_n \rangle /\aaf$ from the
normalized Bar resolution, where $k\langle x_1,\ldots, x_n\rangle$ is the
polynomial ring in $n$ non-commuting indeterminates, and $\aaf$ is a two-sided
ideal with a finite Gr\"obner basis. This result has also been obtained by
Sk\"oldberg \cite{skoed}. In addition to his results we prove the minimality
of this resolution when $\aaf$ is monomial or the Gr\"obner basis consists of
homogeneous polynomials which all have the same degree. In these cases it
follows from our results that the Poincar\'e-Betti series is rational.
In particular, we get the rationality of the Hilbert series if $\aaf$ admits
a quadratic Gr\"obner basis.

In Section \ref{hochcase} we give a projective resolution of $A$ as an
$A\otimes A^{op}$-module, where again
$A=k \langle x_1,\ldots, x_n \rangle /\aaf$. Using this resolution we obtain
the minimal resolution of $A=k[x_1,\ldots, x_n]/\langle f_1,\ldots f_s\rangle$
as an $A\otimes A^{op}$-module, when the initial ideal of
$\langle f_1,\ldots f_s\rangle$ is a complete intersection.
In case $\aaf = \langle f\rangle$ such a construction was first given
by BACH in \cite{bach}.

In the Appendix \ref{appendixa} we derive the normalized Bar and Hochschild
resolution as a sample application of Algebraic Discrete Morse theory and in
Appendix \ref{proof} we give our proof of this theory.

%------------------------------------------------------------------------------
%
% discrete algebraic Morse theory
%
%------------------------------------------------------------------------------
\section{Algebraic discrete Morse theory}\label{secmorse} In this section
we derive an algebraic version of Discrete Morse theory
as developed by Forman (see \cite{forman1}, \cite{forman2}).
Our theory is a generalization of results from \cite{bawel}
and an almost identical theory has been developed independently by
Sk\"oldberg \cite{skoed}. Our applications require a slightly more
general setting than the one covered in \cite{skoed} therefore we
give a detailed exposition of the theory here and provide proofs
in Appendix \ref{proof}.

Let $R$ be a ring and $\CRes_\bullet = (C_i, \partial_i)_{i \geq 0}$ be a
chain complex of free $R$-modules $C_i$. We choose a basis
$X=\bigcup_{i=0}^n X_i$, such that $C_i\simeq \Dirsum_{c\in X_i} R\;c$.
>From now on we write the differentials $\partial_i$ with respect to the
basis $X$ in the following form:

$$
\partial_i: \left\{ \begin{array}{lll} C_i & \to & C_{i-1},\\
       c & \mapsto & \partial_i(c)=\D{\sum_{c'\in X_{i-1}}}[c:c'] \cdot c'.\\
\end{array} \right.
$$

Given the complex $\CRes_\bullet$ and the basis $X$ we construct a directed,
weighted graph $G(\CRes_\bullet)=(V,E)$. The set of vertices $V$
of $G(\CRes_\bullet)$ is the basis $V=X$ and the set $E$ of (weighted) edges
is given by the rule
\begin{eqnarray*}
(c,c',[c:c'])\in E&:\Leftrightarrow& c\in X_i, c'\in X_{i-1}
     \mbox{ and }[c:c']\neq 0.
\end{eqnarray*}
We often omit the weight and write $c \rightarrow c'$ to denote an edge in
$E$. Also by abuse of notation we write $e \in G(\CRes_\bullet)$ to indicate
that $e$ is an edge in $E$.

\begin{Def} \label{morsedefinition}
A subset $\Mmc\subset E$ of the set of edges is called an
acyclic matching, if it satisfies the following three conditions:
\begin{enumerate}
\item (Matching) Each vertex $v\in V$ lies in at most one edge $e\in \Mmc$.
\item (Invertibility) For all edges $(c,c',[c:c'])\in\Mmc$ the weight
  $[c:c']$ lies in the center of $R$ and is a unit in $R$.
\item (Acyclicity) The graph $G_\Mmc(V,E_\Mmc)$ has no directed cycles,
  where $E_\Mmc$ is given by
  \[E_\Mmc:=(E\setminus \Mmc)\cup \left\{\left(c',c,\frac{-1}{[c:c']}\right)
    \mbox{ with }(c,c',[c:c'])\in\Mmc\right\}.\]
\end{enumerate}
\end{Def}

For an acyclic matching $\Mmc$ on the graph $G(\CRes_\bullet) = (V,E)$
we introduce the following notation, which is an adaption of the
notation introduced in \cite{forman1} to our situation.

\begin{enumerate}
\item We call a vertex $c \in V$ critical with respect to $\Mmc$
  if $c$ does not lie in an edge $e\in \Mmc$; we write
  $$X^\Mmc_i:=\{c\in X_i~|~c \mbox{~critical~}\}$$ for the set of all critical
  vertices of homological degree $i$.
\item We write $c'\le c$, if $c\in X_i$, $c'\in X_{i-1}$ and $[c:c']\neq 0$.
\item $\on{Path}(c,c')$ is the set of paths from $c$ to $c'$ in the graph
  $G_\Mmc(\CRes_\bullet)$.
\item The weight $w(p)$ of a path $p=c_1\to\cdots\to c_r\in\on{Path}(c_1,c_r)$
  is given by
  \begin{eqnarray*}
    w(c_1\to\cdots\to c_r)&:=&\prod_{i=1}^{r-1} w(c_i\to c_{i+1}),\\
    w(c\to c')&:=&\left\{\begin{array}{rll}
        -\D{\frac{1}{[c:c']}}&,&c\le c',\\ & & \\
        {\D{[c:c']}}&,&c'\le c.\end{array}\right.
  \end{eqnarray*}
\item We write $\Gamma(c,c')=\D{\sum_{p\in \on{Path}(c,c')}}w(p)$ for the sum
  of weights of all paths from $c$ to $c'$.
\end{enumerate}

\noindent Now we are in position to define a new complex $\CRes^\Mmc_\bullet$
which we call the Morse complex of $\CRes_\bullet$ with respect to $\Mmc$.
The complex $\CRes^\Mmc_\bullet = (C_i^\Mmc,\partial_i^\Mmc)_{i \geq 0}$ is
defined by
\[C_i^\Mmc:=\Dirsum_{c\in X^\Mmc_{i}} R\;c,\]

$$\partial_i^\Mmc: \left\{ \begin{array}{ccc} C_i^\Mmc & \to & C_{i-1}^\Mmc,\\
                c & \mapsto & \D{\sum_{c'\in X^\Mmc_{i-1}}}\Gamma(c,c')c',
\end{array}\right. .$$

\begin{Thm}\label{morse} $\CRes^\Mmc_\bullet$ is a complex of free $R$-modules
and is homotopy equivalent to the complex $\CRes_\bullet$. In particular, for
all $i \geq 0$
\[H_i(\CRes_\bullet) \cong H_i(\CRes^\Mmc_\bullet).\]
The maps defined below induce a chain-homotopy between $\CRes_\bullet$ and
$\CRes^\Mmc_\bullet$.

$$f : \left\{ \begin{array}{lll} \CRes_\bullet & \to     & \CRes^\Mmc_\bullet\\
              c\in X_i             & \mapsto & f(c):=\D{\sum_{c'\in X_i^\Mmc}}\Gamma(c,c')c'
\end{array} \right. $$

$$g : \left\{ \begin{array}{lll} \CRes^\Mmc_\bullet       & \to     & \CRes_\bullet \\
                                c\in X_i^\Mmc & \mapsto & g_i(c):=\D{\sum_{c'\in X_i}}\Gamma(c,c')c'
\end{array} \right. $$

\end{Thm}

\noindent The proof of Theorem \ref{morse} is given in the Appendix
\ref{proof}. Note that if $\CRes_\bullet$ is the cellular chain complex of a
regular CW-complex and $X$ is the set of cells of a regular CW-complex, then
Algebraic Discrete Morse theory is the part Forman's \cite{forman1} Discrete
Morse theory, which describes the impact of a discrete Morse matching on the
cellular chain complex of the CW-complex.

%------------------------------------------------------------------------------
%
% Comm case
%
%------------------------------------------------------------------------------
\section[Resolution of the residue class field]{Resolution of the residue class field
in the commutative case}\label{commcase}

Let $A=S/\aaf$ be the quotient algebra of the commutative polynomial ring
$S= k[x_1, \ldots, x_n]$ in $n$ indeterminates by the ideal $\aaf\subg S$.

The aim of this section is to deduce via Algebraic Discrete Morse theory a new
free resolution of the residue class field $k\iso A/\langle x_1,\ldots ,x_n\rangle$ as an
$A$-module from the normalized Bar resolution. We write $\NBarRes^A_\bullet =
(B_i, \partial_i)_{i \geq}$ for the normalized Bar resolution of $k$ over $A$
(see Appendix \ref{normbar} or \cite{weib}).

>From now on let $\aaf =\langle f_1,\ldots, f_s\rangle\subg S$ be an ideal, such
that the set $\{f_1,\ldots, f_s\}$ is a reduced Gr\"obner basis with respect to a fixed
degree-monomial order `$\prec$' (for example degree-lex or degree-revlex). We assume that
$x_1\succ x_2\succ \ldots \succ x_n$ and we write $\Gc$ for the corresponding set of
standard monomials of degree $\ge 1$.

It is well known, that $\Gc\cup\{1\}$ is a basis of $A$ as $k$-vectorspace.
Thus for any monomial $w$ in $S$ there is a unique representation

\begin{eqnarray} \label{groebnerreduction}
w=a_1+\sum_{v\in \Gc}a_v v, ~ a_1,a_\nu\in k,
\end{eqnarray}

as a linear combination of standard monomials in $A$.

Since we assume that our monomial order is a refinement of the degree order
on monomials it follows that $a_v = 0$ for $|v| > |w|$.
Here we denote with $|v|$ the total degree of the monomial $v$.
In this situation we say that $v$ is reducible to $-\sum_{v\in \Gc}a_v v$.
Note that since we use the normalized Bar resolution the summand $a_1$ can be omitted.

Using the described reduction process we write the
normalized Bar resolution $\NBarRes^A_\bullet = (B_i,\partial_i)$ as:
\begin{eqnarray*}
B_0&:=&A,\\
B_i&:=&\bigoplus_{w_1,\ldots,w_i\in \Gc} A\;[w_1|\ldots|w_i],~i\ge 1
\end{eqnarray*}
with differential
\begin{eqnarray*}
  \lefteqn{\partial_i([w_1|\ldots|w_i])=w_1\,[w_2|\ldots|w_i]}\\
  &&+\D{\sum_{j=1}^{i-1}}(-1)^j \sum_{\nu\in\Gc} a_{j\nu}\:[w_1|\ldots|w_{j-1}|\nu|w_{j+2}\ldots|w_i],
\end{eqnarray*}
for $w_jw_{j+1}=a_{j,1}+\sum_{\nu \in \Gc} a_{j,\nu}\: \nu$, with $a_{j,\nu}\in k$,$\nu \in \Gc $.

The following convention will be convenient.
For a monomial $w \in S$ we set $m(w):=\min\{i\mid x_i\mbox{~divides~}w\}.$
Finally we think of $[w_1|\ldots|w_i]$ as a vector, and we speak of $w_j$ as the entry
in the $j$-th coordinate position.

Now we describe the acyclic matching on the normalized Bar resolution which will be crucial for
the proof of Theorem \ref{k_aufloesung}.
Since all coefficients in the normalized Bar resolutions are $\pm 1$ condition
(Invertibility) of Definition \ref{morsedefinition} is automatically fulfilled.
Thus we only have to take care of the conditions (Matching) and (Acyclicity):

We inductively define acyclic matchings $\Mmc_j$, $j \geq 1$, that are constructed with
respect to the $j$-th coordinate position. We start with the leftmost coordinate position $j=1$.
We set

\[\Mmc_1:=\left\{\begin{array}{c} {[x_{m(w_1)}|w_1'|w_2|\ldots |w_l]} \\
\downarrow \\ {[w_1|w_2|\ldots |w_l]} \end{array} \in
G(\NBarRes_\bullet^A)~\left|~w_1=x_{m(w_1)}w_1'\right.\right\}.\]

The set of critical cells $B^{\Mmc_1}_l$ in homological degree $l\ge 1$ is
given by:

\begin{enumerate}
\item $B^{\Mmc_1}_1:=\Big\{[x_i]~\Big|~ 1\le i\le n\Big\}$, $l=1$
\item $B^{\Mmc_1}_l$ is the set of all $[x_i|w_2|w_3|\ldots|w_l]$, $w_2,\ldots,w_l\in\Gc$, that satisfy
either
  \begin{itemize}
    \item[$\rightarrow$] $i \leq m(w_2)$ and $x_iw_2$ is reducible or
    \item[$\rightarrow$] $i > m(w_2)$.
  \end{itemize}
\end{enumerate}

Assume now $j\ge 2$ and $\Mmc_{j-1}$ is defined. Let $\Bc^{\Mmc_{j-1}}$ be the
set of  critical cells left after applaying $\Mmc_1\cup\ldots\cup\Mmc_{j-1}$.

Let $\Ec_j$ denote the set of edges in $G(\NBarRes_\bullet^A)$ that connect critical cells
in $B^{\Mmc_{j-1}}$.

%On these critical cells and these edges we define for $j \geq 2$ and the $j$-th coordinate:
%\[\Mmc_2:=\Big\{\begin{array}{c} {[x_{i_1}|u_1|u_2|w_3|\ldots|w_k]} \\
%\downarrwo \\ {[x_{i_1}|w_2|w_3|\ldots|w_k]} \end{array} \in \Ec_j~\big|~w_2=u_1u_2\Big\},\]
%where $u_1$ is defined by:
%\[u_1=\max_{\prec}\left\{u\mid w_2\left|\begin{array}{l}[x_{i_1}|u|\ldots]\in B^{\Mmc_1}
%      \\{[x_{i_1}|v|\ldots]}\not\in B^{\Mmc_1}\mbox{, for each divisor $v$ of $u$}\end{array}\right.\right\}.\]
%The set of critical cells is now given by
%\begin{enumerate}
%\item $B^{\Mmc_2}_0:=\Big\{[x_i]~\Big|~1\le i\le n\Big\}$
%\item $B^{\Mmc_2}_k$ consists of elements $[x_i|w_2|w_3|\ldots|w_k]\in B_k^{\Mmc_1}$, such that
%for each divisor $u\mid w_2$ we have $[x_i|u|\ldots]\not\in B_k^{\Mmc_1}$ and one
%  of the following conditions is satisfied:
%  \begin{enumerate}
%    \item  $w_2w_3$ is reducible or
%    \item $w_2w_3=uv\in\Gc$, $[x_i|u|\ldots]\in B_k^{\Mmc_1}$, $u\succ w_2$ and
%      $[x_i|u'|\ldots]\not\in B_k^{\Mmc_1}$ for each divisor $u'\mid u$.
%  \end{enumerate}
%\end{enumerate}
%\lbreak
% In general we define for the $j$-th coordinate:

The following condition on an edge in $\Ec_j$ will define the matching $\Mmc_j$.

\begin{Def}[Matching-Condition]
\label{matchingcondition}
Let $$\begin{array}{c}{[x_{i_1}|w_2|\ldots|w_{j-1}|u_1|u_2|w_{j+1}|\ldots|w_l]}\\ \downarrow\\
{[x_{i_1}|w_2|\ldots |w_{j-1}|w_j|w_{j+1}|\ldots|w_l]}\end{array}$$ be an edge in $\Ec_j$.
In particular $w_j = u_1u_2$. We say that the edge satisfies the matching condition if
$u_1$ is the maximal monomial with respect to `$\prec$' such that
\begin{itemize}
\item[(i)] $u_1$ divides $w_j$.
\item[(ii)] $[x_{i_1}|w_2|\ldots |w_{j-1}|u_1|u_2|w_{j+1}|\ldots|w_l] \in B^{\Mmc_{j-1}}$.
\item[(iii)] $[x_{i_1}|w_2|\ldots |w_{j-1}|v_1|v_2|w_{j+1}|\ldots|w_l]\not\in B^{\Mmc_{j-1}}$ for
for each $v_1 \mid u_1$, $v_1 \neq u_1$ and $v_1v_2 = w_j$.
\end{itemize}
\end{Def}

 \[\Mmc_{j}:=\left\{\begin{array}{c}{[x_{i_1}|w_2|\ldots|w_{j-1}|u_1|u_2|w_{j+1}|\ldots|w_l]}\\
     \downarrow\\{[x_{i_1}|w_2|\ldots |w_{j-1}|w_j|w_{j+1}|\ldots|w_l]}\end{array} \in \Ec_j
~\mbox{satisfying~} \ref{matchingcondition} ~\right\}.\]

We write $\SM(\ini(\aaf))$ for the minimal, monomial generating system
of the initial ideal of $\aaf$ with respect to the chosen monomial order
$\prec$.
The set of critical cells $B^{\Mmc_j}_l$ in homological degree $l\ge 1$ is
given by
\begin{enumerate}
\item $B^{\Mmc_j}_1:=\Big\{[x_i]~\Big|~ 1\le i\le n\Big\}$
\item $B^{\Mmc_j}_2$ consists of elements $[x_i|w_2]$ such
  that either $w_2=x_{i'}$ for some $i'$ and $i>i'$ or $x_iw_2\in\SM(\ini(\aaf))$.
\item $B^{\Mmc_j}_l$ consists of elements
  $[x_i|w_2|\ldots|w_j|\ldots|w_l]\in B_l^{\Mmc_{j-1}}$,
  such that for each divisor $u\mid w_j$ we have
  $[x_i|w_2|\ldots|w_{j-1}|u|\ldots|w_l]\not\in B_l^{\Mmc_{j-1}}$ and one of the following conditions
  is satisfied:
  \begin{itemize}
    \item[$\rightarrow$] $w_jw_{j+1}$ is reducible or
    \item[$\rightarrow$] $w_jw_{j+1}=uv\in\Gc$ and
\begin{itemize}
\item[$\bullet$] $[x_i|w_2|\ldots|w_{j-1}|u|v|w_{j+2}|\ldots|w_l]\in B_l^{\Mmc_{j-1}}$,
\item[$\bullet$] $u\succ w_j$,
\item[$\bullet$] $[x_i|w_2|\ldots|w_{j-1}|u'|v'|w_{j+2}|\ldots|w_l]\not\in B_l^{\Mmc_{j-1}}$ for each
    divisor $u'\mid u$, $u' \neq u$ and $u'v' = w_jw_{j+1}$.
\end{itemize}
\end{itemize}
\end{enumerate}

We finally set $\Mmc:=\bigcup_{j\ge 1}\Mmc_j$ and we write
$\Bc^\Mmc$ for the set of critical cells with respect to $\Mmc$.

\begin{Lem}\label{match}
$\Mmc$ is an acyclic matching.
\end{Lem}
\begin{proof}
We have already seen that since all coefficients are $\pm 1$ the
condition (Invertibility) of Definition \ref{morsedefinition} is
automatic. Property (Matching) is satisfied by definition of
$\Mmc$. Now consider an edge in the matching. Then there exists
a coordinate, where the degree of the monomial decreases by passing
to the higher homological degree cell. Now since we have chosen a
degree-monomial order along any edge in the graph and for any
coordinate the degree of the monomial in this positions
decreases weakly. Since any cycle must contain a matched edge this
shows that there cannot be any directed cycles and (Acyclicity)
is satisfied as well.
\end{proof}
%-------------------------------------------------------------------
%
%Another description of the complex:
%
%-------------------------------------------------------------------
\subsection{An Anick resolution for the commutative polynomial ring}
In this subsection we look closer into the Morse complex corresponding to
the acyclic matching $\Mmc$ from Lemma \ref{match}.
For this we choose the degree-lex order as our fixed monomial order.
We write $\SM(\ini(\aaf))$ for the minimal, monomial generating system
of the initial ideal of $\aaf$ with respect to degree-lex.

In order to describe the critical cells for the chosen term order,
we first define the concept of a minimal fully attached tuple.
Note, that the notation ``fully attached'' was introduced
by Sturmfels (see Example \ref{sturmfels} and \cite{sturm}).

\begin{Def}
A pair $(w_1,w_2)$ is called minimal fully attached, if
$w_1=x_{m(w_1w_2)}$ and $w_1w_2\in \SM(\ini(\aaf))$.

Assume $l>2$. An $l$-tuple $(w_1,\ldots, w_{l-1}, w_l)$ is called minimal
fully attached, if $(w_1,\ldots, w_{l-1})$ is minimal fully attached,
$m(w_1)\le m(w_j)$, for $j=3,\ldots, l$ and one of the following conditions is satisfied:
\begin{enumerate}
\item $w_{l-1}w_l$ is reducible or
\item $w_{l-1}w_l=uv\in\Gc$, with $u\succ w_{l-1}$ and $(w_1,\ldots, w_{l-2},u)$ is a
  minimal fully attached $(l-1)$-tuple,
\end{enumerate}
and $w_l$ is the minimal monomial such that no divisor $w_l'\mid w_l$, $w_l' \neq w_l$
satisfies one of the two conditions above.
\end{Def}

It is easy to see, that the basis of the free modules in the Morse complex $\NBarRes^\Mmc_\bullet$
is given as the set $\Bc$ of words over the alphabet

\begin{eqnarray*}
\Sigma & = & \Big\{ [x_{i_1}|x_{i_2}|\ldots|x_{i_r}]~\Big|~
   1 \leq i_r < i_{r-1} <\cdots <i_1 \leq n~\Big\} \cup \\
       &   & \Big\{ [x_{w_2}|w_2|\ldots|w_l]~|~[x_{w_2}|w_2|\ldots|w_l] \mbox{~minimal fully attached~}\Big\}.
\end{eqnarray*}

that contain none of the words:
\begin{eqnarray}
   [x_{i_1}|\ldots|x_{i_r}][x_{w_2}|w_2|\ldots|w_l], &  x_{w_2} \preceq x_{i_r},\label{rule1}\\[+3mm]
   [x_{i_1}|\ldots|x_{i_r}][x_{j_1}|\ldots|x_{j_s}], & x_{j_1} \preceq x_{i_r},\label{rule2}\\[+3mm]
   [x_{w_2}|w_2|\ldots|w_l][x_{i_1}|\ldots|x_{i_r}], & x_{i_1}\prec x_{w_2},\label{rule3}\\[+3mm]
   [x_{w_2}|w_2|\ldots|w_l][x_{v_2}|v_2|\ldots|v_l], & x_{v_2}\prec x_{w_2}. \label{rule4}
\end{eqnarray}
In order to be able to identify elements of $\Bc$ as basis elements of the Bar resolution we
read in a word from $\Bc$ the sequence of letters `$][$' as `$|$'. If this convention is applied
then any element of $\Bc$ can be read as some $[w_1|\ldots|w_j]$ and corresponds to a basis element in
homological degree $j$. We collect the elements from $\Bc$ which are of homological degree
$j$ in $\Bc_j$ and call an element of $\Bc$ a fully attached tuple.
We claim that there is a bijection between $\Bc^\Mmc$ and $\Bc$ preserving the homological degree.
To see this consider a fully attached tuple $[x_{i_1}|w_2|\ldots|w_i]$.
Then the definition of a fully attached tuple implies, that
either $w_2=x_s$ with $x_s\succ x_{i_1}$ (resp. $i_1>s$) or
$x_{i_1}w_2\in\SM(\ini(\aaf))$. In the first case we cut the tuple to
$[x_{i_1}][x_s|w_3|\ldots|w_i]$. If we continue this process we obtain
\[[x_{i_1}|x_{i_2}|\ldots|x_{i_r}][x_{v_2}|v_2|\ldots|v_s],\]
with $i_1>\ldots>i_r$, $x_{i_r}\prec x_{v_2}$ and
$x_{v_2}v_2\in\SM(\ini(\aaf))$.
This explains the rules (\ref{rule1}) and (\ref{rule2}). Now consider
$[x_{v_2}|v_2|\ldots|v_s]$.
Then the definition of a fully attached tuple implies
that either $v_3=x_j$ with $x_j\succeq x_{v_2}$ or $x_{m(v_3)}\prec x_{v_2}$.
In the first case we cut the tuple to
\[[x_{v_2}|v_2][x_j|v_4|\ldots|v_i],\]
otherwise we consider the monomial $v_4$. Then $v_4$ satisfy the same
conditions as $v_3$, so we cut if necessary to
\[[x_{v_2}|v_2|v_3][x_j|v_5|\ldots|v_i].\]
By construction $[x_{v_2}|v_2|v_3]$ is a minimal fully attached tuple
and the conditions for $v_3$ and $x_j$ explain the rules (\ref{rule3}) and
(\ref{rule4}). If we continue this process we obtain exactly the words
in $\Bc$.

\begin{Rem}\label{poinrem}
Let $\Lc$ be the language over the alphabet
\[\Big\{ [x_{w_2}|w_2|\ldots|w_l]~\Big|~[x_{w_2}|w_2|\ldots|w_l] \mbox{~minimal fully attached~}\Big\},\]
that contain none of the words (\ref{rule4}).
To a letter $[x_{i_1}|x_{i_2}|\ldots|x_{i_r}]\in\Sigma$ with
$1 \leq i_r < i_{r-1} <\cdots <i_1 \leq n$ we associate the symbol
$e_{\{i_r<i_{r-1}<\ldots<i_1\}}$.\\
For $w\in\Bc^{\Mmc}$, such that $w = e_{I_1}\cdots e_{I_s}$, rule (\ref{rule2}) 
shows that this word is, considered as a basis element of $\NBarRes^\Mmc$, equivalent to the
symbol $e_{I_1\cup\ldots\cup I_s}$.\\
To an arbitrary word $w\in\Bc^{\Mmc}$ we first associate the word
\[w_1~e_{I_1}~w_2~e_{I_2}\cdots w_s~e_{I_s}.\]
The rules (\ref{rule1}) and (\ref{rule3}) imply, that the sets $I_i$ are pairwise
disjoint and in a decreasing order. Therefore, as a basis element 
of $\NBarRes^\Mmc$ the word $w$ is equivalent to 
\[e_{I_1\cup\ldots\cup I_s}~w_1w_2\cdots w_s.\]
It follows, that we have a degree-preserving bijection between
$\Bc^{\Mmc}$ and the set
\[\Big\{e_I\wb~\Big|~I\subset \{1,\ldots n\}\mbox{ and }\wb\in \Lc\Big\}.\]
We will use this fact later in order to calculate the multigraded Poincar\'e-Betti series of $k$ over $A$
(see Corollary \ref{poin_cor}).
\end{Rem}
%-------------------------------------------------------------------
%
%description of the differential
%
%-------------------------------------------------------------------

In order to describe the differential, we introduce three reduction rules for fully attached tuples.
These reduction rules will be based on the unique Gr\"obner representation (\ref{groebnerreduction})
which will play role of the basic set of rules:
\[\Rc:=\left\{v_1v_2\stackrel{a_w}{\tol} w\left| \begin{array}{l} v_1,v_2 \in \Gc\\ v_1v_2 \not\in \Gc \end{array} \mbox{~and~} \begin{array}{l}v_1\cdot v_2=a_0+\sum_{w \in \Gc} a_w w,\\
      a_w\in k\end{array}\right.\right\}.\]
Note, that $w\stackrel{0}{\tol} 0\in\Rc$ is allowed (it happens, if one of the generators $f_i$ is a monomial).

\begin{Def}\label{reduction_coeff}
Let $e_1:=[w_1|\ldots|w_{i-1}|w_i|w_{i+1}|w_{i+2}|\ldots| w_l]$ be an $l$-tuple of standard monomials.
\begin{itemize}
\item[{\sf Type I}:]
Assume $[w_1|\ldots|w_i]$ is fully attached. We say $e_1$ can reduced to
$e_2:=[w_1|\ldots|w_{i-1}|v_i|v_{i+1}|w_{i+2}|\ldots| w_l]$, if
\begin{itemize}
\item[(i)] $[w_1|\ldots|w_{i-1}|,v_i]$ is fully attached,
\item[(ii)] $v_iv_{i+1}\in\Gc$,
\item[(iii)] $w_iw_{i+1}\stackrel{a}{\tol}v_iv_{i+1}\in\Rc$, with $a\neq 0$.
\end{itemize}
In this case we write $e_1\stackrel{-a}{\tol}_1 e_2$.
\item[{\sf Type II}:] We say that $e_1$ can
reduced to $e_2:=[w_1|\ldots|w_{i-1}|v|w_{i+2}|\ldots| w_l]$, if
\begin{itemize}
\item[(i)] $w_iw_{i+1}\stackrel{a}{\tol}v\in\Rc$, with $a\neq 0$ and
\item[(ii)] $e_2$ is a fully attached $(l-1)$-tuple.
\end{itemize}
In this case we write $e_1\stackrel{(-1)^ia}{\tol}_2e_2$.
\item[{\sf Type III}:] \label{third} We say that $e_1$ can be reduced to $e_2$
  with coefficient $c:=w_1$ (we write $e_1\stackrel{w_1}{\to}_3 e_2$), if $|w_2|\ge 2$ and
  $e_2:=[x_{m(w_2)}|w_2/x_{m(w_2)}|w_3|\ldots|w_l]$.
\end{itemize}
\end{Def}

Now let $e=[w_1|\ldots| w_l]$ and $f=[v_1|\ldots| v_{l-1}]$ be fully attached $l$- and
$(l-1)$-tuples. We say that $e$ can be reduced to $f$ with coefficient $c$ 
($e\stackrel{c}{\tol}f$) if there exists a sequence $e=e_0, e_1,\ldots,e_{r-1}$ 
and either
\begin{enumerate}
\item an $e_r$, with $e_r=[u|v_1|\ldots|v_{l-1}]=[u|f]$, such that $e_0$ 
  can be reduced to $e_r$ with reductions of Type I and III, i.e.
  \[e_0\stackrel{-a_1}{\tol} e_1\stackrel{-a_2}{\tol}e_2\stackrel{-a_3}{\tol}\ldots\stackrel{-a_r}{\tol}e_r;\]
  in this case we set $c:=\left((-1)^r\prod_{i=1}^ra_i\right)\,u$, or
\item an $e_r$, such that $e_0$ can be reduced to $e_r$ with reductions of Type I and III and
  $e_r$ can be reduced to $f$ with reduction of Type II, i.e.
  \[e_0\stackrel{-a_1}{\tol} e_1\stackrel{-a_2}{\tol}e_2\stackrel{-a_3}{\tol}\ldots\stackrel{-a_r}{\tol}e_r\stackrel{(-1)^jb}{\tol} f;\]
  in this case, we set $c:=(-1)^{r+j}\cdot b\cdot\prod_{i=1}^ra_i$.
\end{enumerate}

There may be several possible reduction sequences leading from $e$ to $f$ and the reduction coefficient
may depend chosen sequence. Therefore, we define {\em the} reduction
coefficient $[e:f]$ to be the sum over all possible sequences. If there exists no sequence, we set
$[e:f]:=0$.

%-------------------------------------------------------------------
%
%end of description of the complex.
%
%-------------------------------------------------------------------
The complex $\FRes_\bullet$ is then given by
\[F_j:=\bigoplus_{e\in \Bc_j}A\,e,\]
\begin{eqnarray*}
\partial:F_i&\to&F_{i-1}\\
e&\mapsto&\sum_{f\in \Bc_{i-1}}[e:f]\, f.
\end{eqnarray*}

Now we have:

\begin{Thm}\label{k_aufloesung}
$\FRes_\bullet = (F_\bullet, \partial)$ is an $A$-free resolution of the residue class field $k$, which is minimal
if and only if no reduction of Type II is possible.
\end{Thm}
\begin{proof} The fully attached tuples are exactly the critical cells. The reduction rules
describe the Morse differential: As seen before, we have
\[\partial^\Mmc([w_1,\ldots, w_l]):=w_1[w_2,\ldots, w_l]+\sum_{i=1}^{l-1}(-1)^i[w_1,\ldots, w_iw_{i+1}, \ldots, w_l]\]
If $[w_2,\ldots, w_l]\not\in\Bc$, we have $[w_2,\ldots, w_l]=\partial([x_{i_2}, w_2',w_3\ldots,
w_l])$, which is described by the reduction of Type III.

For $[w_1,\ldots, w_iw_{i+1}, \ldots, w_l]$ we have to distinguish three cases:\\
\begin{itemize}
\item[(Case 1)] $[w_1,\ldots, v_{ij}, \ldots, w_l]$ is critical. Then we have $w_{i-1}v_{ij}$,
$v_{ij}w_{i+2}$ reducible and $w_{i-1}u_1$ $v_{ij}u_2\in\Gc$ for all divisors $u_1$ of
$v_{ij}$ and $u_2$ of $w_{i+2}$. This situation is described by the reduction of Type II.
\item[(Case 2)] $[w_1,\ldots, v_{ij}, \ldots, w_l]$ is matched by a higher degree cell. Then we
have: $w_{i-1}u_1$ reducible for $v_{ij}=u_1u_2$ and for all divisors $u'$ of $u_1$ the
monomial $w_{i-1}u'$ lies in $\Gc$. Then we have $$[w_1,\ldots, v_{ij}, \ldots, w_l]=
(-1)^{i+1}[w_1,\ldots, w_{i-1},u_1,u_2,w_{i+2}, \ldots, w_l],$$ which is a reduction
of Type I.
\item[(Case 3)] $[w_1,\ldots, v_{ij}, \ldots, w_l]$ is matched by an lower degree cell. In this case
we have $[w_1,\ldots, v_{ij}, \ldots, w_l]=0$.
\end{itemize}
The coefficients of the reductions are exactly the coefficients of the Morse differential. Hence
the Morse-differential induces a sequence of reductions of Type I and III with either a reduction of
Type II, or the map $e_r=[v_1|\ldots|v_l]\stackrel{v_1}{\tol} [v_2\ldots,v_l]$ at the end, which
gives our definition of the reduction coefficient.
\end{proof}
\begin{Rem}
In Section \ref{noncommcase}, we will see, that in the non-com\-mutative case our matching on the
normalized Bar resolution gives the Anick resolution (for the definition, see \cite{Anick}).
Therefore one can understand the resolution $F_\bullet$ as a generalization of the Anick resolution to the
commutative polynomial ring.
\end{Rem}
If $A$ is endowed with the natural multigrading $\deg(x_i)=e_i\in\N^n$ the multigraded Poincar\'e-Betti series
of $k$ over $A$ is defined to be
\[\poinr{k}{A}:=\sum_{i\ge 0\atop \alpha\in\N^n}\dim_k(\Tor^A_i(k,k)_\alpha)\,\ul{x}^\alpha\,t^i.\]
Remark \ref{poinrem} implies:

\begin{Cor}\label{poin_cor}
The Poincar\'e-Betti series of $A$ satisfies
\[\poinr{k}{A}(x,t)\le \prod_{i=1}^n(1+x_i\,t)\;F(x,t),\]
where $F(x,t):=\sum_{w\in\Lc}w\,t^{|w|}$ counts the words $w\in\Lc$. Here $w$ is treated as the monomial in $x_1,\ldots,x_n$ and $|w|$ denotes the length of $w$.\\
The inequality is an inequality between the
coefficients of the power series expansion. \qed
\end{Cor}
%-------------------------------------------------------------------
%
%two special cases
%
%-------------------------------------------------------------------
\subsection{Two special cases}
First we consider a subclass of the class of Koszul algebras.
It is well known, that $A=S/\aaf$ is Koszul if $\aaf$ has a
quadratic Gr\"obner basis. It is easy to see, that in this case the minimal fully attached tuples have
the following form: $[x_{i_1}|x_{i_2}|\ldots|x_{i_r}]$. Therefore a reduction of Type II is not possible
and we get:

\begin{Cor}\label{koszul_minimal}
If $A=S/\aaf$ and $\aaf$ admits a quadratic Gr\"obner basis, then
the resolution $F_\bullet$ is minimal.\qed
\end{Cor}
To get an explicit form of the multigraded Poincar\'e-Betti series in this case one only has to
calculate the word-counting function $F(x,t)$ of the language $\Lc$.
In this case the multigraded Poincar\'e-Betti series coincides with the multigraded Poincar\'e-Betti
series of $S/\ini(\aaf)$. Since the Poincar\'e-Betti series of monomial rings are studied by us in a larger
context in \cite{jol} we do not give the explicit form here.

\medskip

The second case, we would like to discuss, is the following:\\
Let $\aaf=\langle f_1,\ldots, f_s\rangle\subg S$ be an ideal, such
that $f_1,\ldots, f_s$ is a reduced Gr\"obner basis with respect
to the degree-lex order and such that the initial ideal
$\ini(\aaf)$ is a complete intersection. Assume
$f_j=m_j+\sum_{\alpha\in\mathbb{N}^n}f_{j\alpha}x^\alpha$ with
leading monomial $m_j$. Since $\ini(\aaf)$ is a complete
intersection, there exist exactly $s$ minimal fully attached
tuples, namely
$t_i:=\left[x_m(m_i)\left|\frac{m_i}{x_m(m_i)}\right.\right]$, for
$i=1,\ldots, s$ and $m_i\in \SM(\ini(\aaf))$. The rule
(\ref{rule4}) implies $t_it_j\in\Bc$, iff $m(m_i)\ge m(m_j)$. It
follows from Remark \ref{poinrem}, that the set of fully attached
$i$-tuples is in bijection with the set
\[\Bc_i:=\left\{e_{i_r}\ldots e_{i_1}t_{j_1}^{(l_1)}\ldots t_{j_q}^{(l_q)}\left|
  \begin{array}{l}1\le i_1<\ldots<i_r\le n\\1\le j_1<\ldots<j_q\le s\\l_1,\ldots, l_q\in\N\mbox{ and }i=r+2\sum_{t=1}^ql_t
  \end{array}\right.\right\}.\]
%
%In order to describe the differential, we further assume, that for all $f_{i\alpha}\neq 0$ and all $j\neq i$ we have
%\[m_j\not|\; x_r\frac{x^{\alpha}}{x_{m(x^{\alpha})}}\mbox{, for all $r=1,\ldots, n$}.\]
For $f_j=m_j+\sum_{\alpha\in\N^n}f_{j\alpha}x^\alpha$ we define
\[T_p(f_j):=\sum_{\alpha\in\N^n\atop p=\max(\supp(\alpha))}f_{j\alpha}\;\frac{x^{\alpha}}{x_p}.\]
We have the following theorem:
%-------------------------------------------------------------------
%
% theorem complete intersection
%
%-------------------------------------------------------------------
\begin{Thm}\label{bachkomm} Let $\aaf=\langle f_1,\ldots, f_s\rangle\subg S$ be an ideal, such
that $f_1,\ldots, f_s$ is a reduced Gr\"obner basis with respect
to the degree-lex order and such that the initial ideal
$\ini(\aaf)$ is a complete intersection and $A:=S/\aaf$ be the quotient algebra.\\
Then the following complex is a minimal $A$-free resolution of the residue class field $k$
and carries the structure of a differential graded algebra:
\[F_i:=\bigoplus_{\begin{array}{c}1\le i_1<\ldots<i_r\le n\\1\le j_1<\ldots<j_q\le s\\
    l_1,\ldots, l_q\in\mathbb{N}\\i=r+2\sum_{j=1}^ql_j\end{array}}A\; e_{i_r}\ldots e_{i_1}t_{j_1}^{(l_1)}\ldots t_{j_q}^{(l_q)}\]
\begin{eqnarray*}
e_{i_r}\ldots e_{i_1}&\stackrel{\partial}{\mapsto}&\sum_{m=1}^r (-1)^{\#\{i_j>i_m\}}x_{i_m}\,e_{i_r}\ldots \widehat{e_{i_m}}\ldots e_{i_1}\\
t_{j_1}^{(l_1)}\ldots t_{j_q}^{(l_q)}&\stackrel{\partial}{\mapsto}&
   \sum_{m=1}^s\sum_{p=1}^n\,T_p(f_{j_m})\;e_pt_{j_1}^{(l_1)}\ldots, t_{j_m}^{(l_{j_m}-1)},\ldots, t_{j_q}^{(l_q)},
\end{eqnarray*}
where $t_{i_j}^{(0)}:=1$, $e_ie_j=-e_je_i$ and $e_ie_i=0$. The differential is given by
\begin{eqnarray*}
\partial(e_{i_r}\ldots e_{i_1}t_{j_1}^{(l_1)}\ldots t_{j_s}^{(l_s)})&=&\partial(e_{i_r}\ldots e_{i_1})t_{j_1}^{(l_1)}\ldots t_{j_s}^{(l_s)}\\
&&(-1)^r\,e_{i_r}\ldots e_{i_1}\partial(t_{j_1}^{(l_1)}\ldots t_{j_s}^{(l_s)}).
\end{eqnarray*}
In particular we have
\[\poinr{k}{A}(x,t)= \frac{\D\prod_{i=1}^n(1+x_i\,t)}{\D \prod_{i=1}^k(1-m_i\,t^2)}.\]
\end{Thm}
\begin{proof}
We only have to calculate the differential: Let $[w_1,\ldots,
w_l]$ be a fully attached tuple, such that $w_j$ is either a
variable or a minimal fully attached tuple.\\
First assume that $w_j$ is a variable, i.e. $w_j=x_{r_j}$. We
prove that $w_iw_j$ can be permuted to $w_jw_i$ for all $i\neq
j$. If $w_i$ is a variable, say $w_i=x_{j_i}$, we have by
(\ref{rule4}) $j_i>r_j$ it follows $|x_{j_i}|x_{j_r}|\to
|x_{j_i}x_{r_j}|\to |x_{r_j}|x_{j_i}|$. If $w_i$ is a minimal
fully attached tuple, i.e.
$w_i=\left|x_{m(m_i)}\left|\frac{m_i}{x_{m(m_i)}}\right.\right|$,
we have
\begin{eqnarray*}
\left|x_{m(m_i)}\left|\left.\frac{m_i}{x_{m(m_i)}}\right|x_{r_j}\right.\right|&\to&\left|x_{m(m_i)}
  \left|x_{r_j}\frac{m_i}{x_{m(m_i)}}\right.\right|\to\left|x_{m(m_i)}\left|x_{r_j}\left|\frac{m_i}{x_{m(m_i)}}\right.\right.\right|\\
&\to&\left|x_{r_j}x_{m(m_i)}\left|\frac{m_i}{x_{m(m_i)}}\right.\right|\to
\left|x_{r_j}\left|x_{m(m_i)}\left|\frac{m_i}{x_{m(m_i)}}\right.\right.\right|
\end{eqnarray*}
In the first case we have a reduction with coefficient $-1$ and in the second case
 with coefficient $+1$.
Therefore it is enough to consider the number of $w_i's$, $i<j$,
which are variables. It follows, that $w_j$ can be permuted to the
left with coefficient $(-1)^{\#\{w_i~|~w_i~\mbox{variable and
}~w_i<_{\lex}x_{r_j}\}}$.\\
Now let $w_j$ be a minimal fully attached tuple, i.e.
$w_j=\left[x_{m(m_j)}\left|\frac{m_j}{x_{m(m_j)}}\right.\right]$.
Then we have
\[\left[x_{m(m_j)}\left|\frac{m_j}{x_{m(m_j)}}\right.\right]\to -\sum_{\alpha}f_{j\alpha}[x^{\alpha}]
  \to \sum_{\alpha}f_{j\alpha}\left[x_{\alpha}\left|\frac{x^{\alpha}}{x_{\alpha}}\right.\right],\]
where $x_\alpha:=x_{m(x^\alpha)}$.
Since $\left[\frac{x^{\alpha}}{x_{\alpha}}\right]$ is matched with
$\left[x_{\beta}\left|\frac{x^{\alpha}}{x_{\beta}x_{\alpha}}\right.\right]$
(where $x_\beta=x_{m(x^\beta)}$ with $x^\beta:=\frac{x^{\alpha}}{x_{\alpha}}$) 
the exponent $\alpha$ decreases
successively up to the element $[x_p]$, with $p=\max(\supp(\alpha))$. Therefore we get
\begin{equation}
\left[x_{m(m_j)}\left|\frac{m_j}{x_{m(m_j)}}\right.\right]\to \sum_{p=1}^nT_p(f_j)e_p.
\label{complete_diff}
\end{equation}
We now consider the tuple $[w_1,\ldots, w_l]$. With the same
argument as before one can check, that the minimal fully attached tuple $w_j$ can be permuted with
coefficient $+1$ to the right. After a chain of reductions we reach the tuple
$[w_j,w_1,\ldots,w_{j-1},w_{j+1},\ldots,w_l]$. Applying Equation (\ref{complete_diff})
we get
\[[w_1,\ldots, w_l]\to\sum_{p=1}^nT_p(f_j)[x_p,w_1,\ldots, \widehat{w_j},\ldots,w_l].\]
In order to reach a fully attached tuple we have to permute the variable $x_p$ to the correct position.
This permutation yields a coefficient $(-1)^{\#\{w_i~|~w_i~\mbox{variable and}~w_i<_{\lex}x_p\}}$.\\
The bijection between the elements $e_{i_r}\ldots e_{i_1}t_{j_1}^{(l_1)}\ldots t_{j_q}^{(l_q)}$ 
and the fully attached tuples finally implies the coefficient
$$\Big((-1)^{\#\{w_i~|~w_i~\mbox{variable and}~w_i<_{\lex}x_p\}}\Big)^2
(-1)^r=(-1)^r.$$
Therefore our differential has the desired form
\begin{eqnarray*}
\lefteqn{\partial(e_{i_r}\ldots e_{i_1}t_{j_1}^{(l_1)}\ldots t_{j_q}^{(l_q)})}\\
&=&\sum_{m=1}^r (-1)^{\#\{i_j>i_m\}}x_{i_m}\,e_{i_r}\ldots \widehat{e_{i_m}}\ldots e_{i_1}t_{j_1}^{(l_1)}\ldots t_{j_q}^{(l_q)}\\
&&+\sum_{m=1}^q\sum_{p=1\atop p\neq i_1,\ldots,i_r}^n(-1)^r\:T_p(f_{j_m})\;
e_{i_r}\cdots e_{i_1}e_pt_{j_1}^{(l_1)}\ldots, t_{j_m}^{(l_{j_m}-1)},\ldots, t_{j_q}^{(l_q)}.
\end{eqnarray*}
It is easy to see, that this are all possible reductions.
\end{proof}

If $\ini(\aaf)=\aaf$ then the preceding result about the Poincar\'e-Betti series 
cab be found in \cite{Gull}.
%------------------------------------------------------------------------------
%
% Non-comm case
%
%------------------------------------------------------------------------------
\section[Resolution of the residue class field]{Resolution of the residue class field in
the non-commutative case}\label{noncommcase}

In this section we study the same situation as in Section \ref{commcase}
over the polynomial ring in $n$ non-commuting indeterminates.
In this case the acyclic matching on the normalized Bar resolution is
slightly different to the acyclic matching in Section \ref{commcase} and the
resulting Morse complex will be isomorphic to the Anick resolution. These
results were independently obtained by Sk\"oldberg \cite{skoed}. In addition to 
Sk\"oldberg's results we prove minimality of this resolution in special
cases, which give information about the Poincar\'e-Betti series, and we give
an explicit description of the complex if the two-sided ideal $\aaf$
admits a (finite) quadratic Gr\"obner basis, which proves a conjecture by
Sturmfels \cite{sturm}.

\medskip

Let $A=k \langle x_1, \ldots, x_n \rangle/\aaf$ be the quotient algebra
of the polynomial ring in $n$ non-commuting indeterminates by a
two-sided ideal $$\aaf\subg k \langle x_1, \ldots, x_n \rangle.$$
As before, we assume, that $\aaf=\langle f_1,\ldots f_s\rangle$, such that
$\{f_1,\ldots, f_s\}$ is a finite reduced Gr\"obner basis with respect to a
fixed degree-monomial order $\prec$. For an introduction to the
theory of Gr\"obner basis in the non-commutative case see \cite{groebner}.

Again we have for the product of any two standard monomials a unique
(Gr\"obner-) representation of the form:
\[w\cdot v:=\sum_i a_i\, w_i\mbox{ with $a_i\in k$, $w_i\in\Gc$
  and $|w\cdot v|\ge |w_i|$ for all $i$},\]
where $\Gc$ is the corresponding set of standard monomials of degree
$\ge 1$ and $|m|$ is the total degree of the monomial $m$.

The acyclic matching on the normalized Bar resolution is defined as follows:
As in the commutative case, we define $\Mmc_j$ by induction on the coordinate $1\leq j \leq n$: 
For $j=1$ we set
\[\Mmc_1:=\left\{\begin{array}{c}{[x_i|w_1'|w_2|\ldots|w_l]}\\\downarrow\\
    {[w_1|\ldots|w_l]}\end{array}~\in G(\NBarRes_\bullet^A)
  \left|~ w_1=x_iw_1'\right.\right\}.\]
The critical cells with respect to $\Mmc_1$ are given by
\begin{enumerate}
\item $B^{\Mmc_1}_1:=\Big\{[x_i]~\Big|~ 1\le i\le n\Big\}$, $l=1$
\item $B^{\Mmc_1}_l$ is the set of all $[x_i|w_2|w_3|\ldots|w_l]$, $w_2, w_3,\ldots,w_l\in\Gc$, 
  that satisfy 
  \begin{itemize}
    \item[$\rightarrow$]$x_iw_2$ is reducible
  \end{itemize}
\end{enumerate}

Assume now $j\ge 2$ and $\Mmc_{j-1}$ is defined. Let $\Bc^{\Mmc_{j-1}}$ be the
set of  critical cells left after applaying $\Mmc_1\cup\ldots\cup\Mmc_{j-1}$.

Let $\Ec_j$ denote the set of edges in $G(\NBarRes_\bullet^A)$ that connect critical cells
in $B^{\Mmc_{j-1}}$.

The following condition on an edge in $\Ec_j$ will define the matching $\Mmc_j$.

\begin{Def}[Matching-Condition]
\label{matchingcondition_noncom}
Let $$\begin{array}{c}{[x_{i_1}|w_2|\ldots|w_{j-1}|u_1|u_2|w_{j+1}|\ldots|w_l]}\\ \downarrow\\
{[x_{i_1}|w_2|\ldots |w_{j-1}|w_j|w_{j+1}|\ldots|w_l]}\end{array}$$ be an edge in $\Ec_j$.
In particular $w_j = u_1u_2$. We say that the edge satisfies the matching condition if
\begin{itemize}
\item[(i)] $u_1$ is a prefix of $w_j$.
\item[(ii)] $[x_{i_1}|w_2|\ldots |w_{j-1}|u_1|u_2|w_{j+1}|\ldots|w_l] \in B^{\Mmc_{j-1}}$.
\item[(iii)] $[x_{i_1}|w_2|\ldots |w_{j-1}|v_1|v_2|w_{j+1}|\ldots|w_l]\not\in B^{\Mmc_{j-1}}$ for
for each prefix $v_1$ of $u_1$ and $v_1v_2 = w_j$.
\end{itemize}
\end{Def}
 \[\Mmc_{j}:=\left\{\begin{array}{c}{[x_{i_1}|w_2|\ldots|w_{j-1}|u_1|u_2|w_{j+1}|\ldots|w_l]}\\
     \downarrow\\{[x_{i_1}|w_2|\ldots |w_{j-1}|w_j|w_{j+1}|\ldots|w_l]}\end{array} \in \Ec_j
~\mbox{satisfying~} \ref{matchingcondition_noncom} ~\right\}.\]

The set of critical cells $\Bc_l^{\Mmc_j}$ in homological degree $l\ge 1$ is given by
\begin{enumerate}
\item $B^{\Mmc_j}_1:=\Big\{[x_i]~\Big|~ 1\le i\le n\Big\}$
\item $\D B^{\Mmc_j}_2$ consists of elements $[x_{i_1}|w_2]$, with
  $x_{i_1}w_2\in \SM(\ini(\aaf))$.
\item $B^{\Mmc_j}_l$ consists of elements $[x_{i_1}|w_2|w_3|\ldots|w_l]\in\Bc^{\Mmc_{j-1}}_l$,
  such that for each prefix $u$ of $w_j$ we have 
  $[x_{i_1}|w_2|\ldots|w_{j-1}|u|\ldots w_l]\not\in\Bc^{\Mmc_{j-1}}_l$
  and $w_jw_{j+1}$ is reducible.
\end{enumerate}

We finally set $\Mmc:=\bigcup_{j\ge 1}\Mmc_j$ and we write $\Bc^\Mmc$ for the set of critical
cells with respect to $\Mmc$.

With the same proof as in Section \ref{commcase} we get
\begin{Lem}\label{match2} $\mathcal{M}$ defines an acyclic matching.\qed
\end{Lem}
%-------------------------------------------------------------------
%
%Another description of the complex:
%
%-------------------------------------------------------------------
\subsection{The Anick resolution}
As in the commutative case we give a second description of the Morse complex
with respect to the acyclic matching from Lemma \ref{match2}.
In this case this description shows that it is isomorphic to the 
Anick resolution \cite{Anick}.

\begin{Def} Let $m_{i_1},\ldots, m_{i_{l-1}}\in \SM(\ini(\aaf))$ be monomials,
such that for $j=1,\ldots l-1$ we have $m_{i_j}=u_{i_j}v_{i_j}w_{i_j}$ with
$u_{i_{j+1}}=w_{i_j}$ and $|u_{i_1}|=1$. Then we call the $l$-tuple
$$[u_{i_1},v_{i_1}w_{i_1},v_{i_2}w_{i_2},\ldots, v_{i_{l-1}}w_{i_{l-1}}]$$
fully attached, if for all $1\le i\le l-2$ and each prefix $u$ of
$v_{i_{j+1}}w_{i_{j+1}}$ the monomial $v_{i_j}w_{i_j}u$ lies in $\Gc$.
We write $\Bc_j:=\{[w_1,\ldots,w_j]\}$ for the set of fully attached $j$-tuples
($j\ge 2$) and $\Bc_1:=\{[x_1],\ldots, [x_n]\}$.
\end{Def}

We define the reduction types (Type I, Type II and Type III)
and the reduction coefficient $[e:f]$ for two fully attached tuples $e,f$ 
in a similar way as in the commutative case (see Definition \ref{reduction_coeff}).
Now we are able to define the following complex:
\[F_j:=\Dirsum_{e\in \Bc_j}A\,e,\]
\begin{eqnarray*}
\partial:F_i&\to&F_{i-1}\\
e&\mapsto&\sum_{f\in \Bc_{i-1}}[e:f]\, f.
\end{eqnarray*}

Note, that the basis elements of $F_j$ are exactly the basis elements in 
the Anick resolution (see \cite{Anick}), therefore the complex $F_\bullet$ is 
isomorphic to the Anick resolution.
Again we have:
\begin{Thm}\label{k_res_noncomm}
$(F_\bullet, \partial)$ is an $A$-free resolution of the residue class field $k$ over $A$.
If no reduction of Type II is possible, the resolution $(F_\bullet,\partial)$ is minimal.
\end{Thm}
\begin{proof} The fully attached tuples are exactly the critical cells. 
The rest is analogous to the commutative case.
\end{proof}
If one applies Theorem \ref{k_res_noncomm} to the ideal $\langle x_ix_j-x_jx_i,\aaf\rangle$, 
one reaches the commutative case. But in general the Morse complex with respect to the 
acyclic matching from Lemma \ref{match2} is much larger (with respect to the rank) 
than the Morse complex of the acyclic matching, developed in Section \ref{commcase} 
for commutative polynomial rings.\\

Since only by reductions of Type II coefficient $[e:f]\in k$ can 
enter the resolution, we have:
\begin{Prop}\label{minbed}The following conditions are equivalent:
\begin{enumerate}
\item $(F_\bullet,\partial)$ is not minimal.
\item There exist standard monomials $w_1,\ldots w_4$, and
  $m_{i_1},m_{i_2},m_{i_3}\in \SM(\ini(\aaf))$, such that
  $w_1w_2=u_1m_{i_1}$ ,$w_2w_3=u_2m_{i_2}$, $w_1w_4=u_1'm_{i_3}$ with
  $u_1,u_1'$ suffix of $w_1$, $u_2$ suffix of $w_2$  and $w_2w_3\to w_4\in \Rc$
\end{enumerate}
\end{Prop}
\begin{proof}
$(F_\bullet,\partial)$ is minimal iff no reduction of Type II is possible,
which is equivalent to the second condition.
\end{proof}
\begin{Cor}\label{red}
In the following two cases the resolution $(F_\bullet,\partial)$ is a minimal
$A$-free resolution of $k$ and independent of the characteristic of $k$.
\begin{enumerate}
\item $\aaf$ admits a monomial Gr\"obner basis.
\item The Gr\"obner basis of $\aaf$ consists of homogeneous polynomials, all
  of the same degree.
\end{enumerate}
\end{Cor}
\begin{proof}
If the Gr\"obner basis consists of monomials, the situation of Proposition \ref{minbed} is not possible.
In the other case there exists a constant $l$, such that for all $w\to v\in\mathcal{R}$ we have
$|w|=|v|=l$. Assume there exist standard monomials $w_1,\ldots, w_4$, and monomials
$m_{i_1},m_{i_2},m_{i_3}\in \SM(\ini(\aaf))$, such that $w_1w_2=u_1m_{i_1}$ ,$w_2w_3=u_2m_{i_2}$,
$w_1w_4=u_1'm_{i_3}$ with $u_1,u_1'$ suffix of $w_1$,  $u_2$ suffix of $w_2$ and
$w_2w_3\to w_4\in \mathcal{R}$. Then we get $|w_i|<l$ for $i=2,3,4$. On the other hand we have
$w_2w_3\to w_4\in \mathcal{R}$ and therefore $|w_4|=l$. This is a contradiction.
\end{proof}
%-------------------------------------------------------------------
%
%Poincar\'e series of $k$
%
%-------------------------------------------------------------------
\subsection{The Poincar\'e-Betti series of $k$}
In this subsection we draw some corollaries on the Poincar\'e-Betti series of
$k$.\\

Recall the definition of a fully attached $l$-tuple: There exist leading
monomials $m_{i_1},\ldots,m_{i_{l-1}}\in \SM(\ini(\aaf))$, such that for all
$j=1,\ldots l-1$ there exists a monomial
$u_{i_j},v_{i_j},w_{i_j}\in\mathcal{G}$ with $m_{i_j}=u_{i_j}v_{i_j}w_{i_j}$
and $u_{i_{j+1}}=w_{i_j}$. It follows that the fully attached $l$-tuples are
in one to one correspondence with $l-1$ chains of monomials
$(m_{i_1},\ldots, m_{i_{l-1}})$ with the condition before. We write again
$\Bc$ for the set of all those chains. Now consider the subset
$$E:=\Big\{(m_{i_1}, \ldots, m_{i_l})\in\Bc~\Big|~ m_{i_1},\ldots, m_{i_l}
  \mbox{ pairwise different }\Big\} \subset \Bc.$$
Since we consider only finite Gr\"obner basis it is clear, that $E$ is
finite.
We construct a DFA (deterministic finite automaton, see for example
\cite{automat}) over the alphabet $E$, which accepts $\Bc$. 
For each letter $f\in E$ we define a
state $f$. Each state $f$ is a final state.
Let $S$ be the initial state, and $Q$ be a error state.
>From the state $S$ we go to state $f$ if we read the letter $f\in E$.
Let $f_1=(m_{i_1}, \ldots, m_{i_l}),f_2=(m'_{j_1}, \ldots, m'_{j_{l'}})\in E$ be
two chains of monomials with corresponding fully attached tuples
$(w_{i_1},\ldots, w_{i_{l+1}})$ and $(w'_{j_1},\ldots, w'_{j_{l'+1}})$.
Then we have $(f_1,f_2)\in\Bc$, iff there exists a monomial
$n\in \SM(\ini(\aaf))$ with $n=uw'_{j_1}$ and $u$ suffix of $w_{i_{l+1}}$.
In this case we change by reading $f_2$ from state $f_1$ to $f_2$.
If such a monomial does not exist we change  by reading $f_2$ from state
$f_1$ to the error state $Q$. The language of this DFA is exactly the set
$\Bc$. This proves that the basis of our resolution $F_\bullet$
is a regular language. Since the word-counting function of a regular language
is always a rational function (see \cite{automat}) we get in particular the
following theorem:
\begin{Thm}\label{pointhm_noncomm} For the Poincar\'e-Betti series of $k$ over
the ring $A$ we have \[\poinr{k}{A}(x,t)\le F(x,t),\]
where $F(x,t)$ is a rational function. Equality holds iff $F_\bullet$ is
minimal.\qed
\end{Thm}
\begin{Cor}
For the following two cases the Poincar\'e-Betti series of $k$ over the ring
$A$ is a rational function:
\begin{enumerate}
\item $\aaf$ admits a Gr\"obner basis consisting of monomials.
\item The Gr\"obner basis of $\aaf$ consists of homogeneous polynomials, all
  of the same degree.
\end{enumerate}
\end{Cor}
\begin{proof}
The result is a direct consequence of the theorem \ref{pointhm_noncomm} and
Corollary \ref{red}.
\end{proof}
\begin{Cor} If $\aaf$ has a quadratic Gr\"obner basis, then $F_\bullet$ is an
  $A$-free minimal linear resolution. Hence $A$ is Koszul and its Hilbert- and
  Poincar\'e-Betti  series are rational functions.\qed
\end{Cor}

%-------------------------------------------------------------------
%
%Examples
%
%-------------------------------------------------------------------
\subsection{Examples} We finally want to give some examples of the Morse complex and 
we verify a conjecture by Sturmfels:

\begin{Ex}[Conjecture of B. Sturmfels (see \cite{sturm})]\label{sturmfels}
Let $\Lambda$ be a graded subsemigroup of $\N^d$ with $n$ generators.
We write its semigroup algebra over a field $k$ as a quotient of the free associative
algebra $$ k \langle y_1,y_2,\ldots, y_n\rangle / J_\Lambda\quad = \quad k[ \Lambda ].$$
Suppose that the two-sided ideal $J_\Lambda$ possesses a quadratic Gr\"obner basis  $\Gc$.
The elements in the non-commutative Gr\"obner basis $\Gc$
are quadratic reduction relations of the form
$\,y_i y_j  \rightarrow y_{i'} y_{j'}$. If $w$ and $w'$ are words in  $y_1,\ldots,y_n$ then we write
$ \, w {\buildrel j \over \longrightarrow} w' \,$ if there exists
a reduction sequence of length $j$ from $w$ to $w'$.
A word $\,w =  y_{i_1} y_{i_2} \cdots y_{i_r} \,$
is called {\it fully attached} if every quadratic subword
$\,y_{i_j} y_{i_{j+1}} \,$ can be reduced with respect to $\Gc$.
Let ${\bf F}_r$ be the free $k[\Lambda]$-module with basis
$\,\{ \,E_w \,: \, w \,$ fully attached word of length $r \}$.
Let ${\bf F} =  \bigoplus_{r \geq 0} {\bf F}_r$ and define
a differential $\partial $ on ${\bf F}$ as follows:
$$\,\partial : {\bf F}_r \rightarrow {\bf F}_{r-1} \, ,\quad
E_w  \quad \mapsto \quad
\sum (-1)^j \,x_i \,E_{w'}\,$$
where the sum is over all fully attached words
$w'$ of length $r-1$ such that
$ \, w \,{\buildrel j \over \longrightarrow}  \,x_i\, w'$
for some $i,j$. Note that this sum includes the
trivial reduction  $\,w {\buildrel 0 \over \longrightarrow}  w $.
Then Theorem \ref{k_res_noncomm} together with Proposition
\ref{minbed} implies that $({\bf F}, \partial)$ is a minimal free resolution of $k$ over $k[\Lambda]$.\qed
\vskip .3cm
\noindent {\em Example 1} {\sl (The twisted cubic curve): }
The Gr\"obner basis  consists of nine binomials:
$$ {\mathcal G} \,= \, \left\{\begin{array}{l}
ac \rightarrow bb ,\,
ca \rightarrow bb ,\,
ad \rightarrow cb ,\,
da \rightarrow cb ,\,
bd \rightarrow cc ,\\
db \rightarrow cc ,\,
ba \rightarrow ab ,\,
bc \rightarrow cb ,
dc \rightarrow cd \end{array}\right\} $$
The minimal free resolution $({\bf F},\partial)$ has the format
$$ \cdots \cdots
\,\, {\buildrel \partial \over \longrightarrow}  \,\,
k[\Lambda]^{72} \,\, {\buildrel \partial \over \longrightarrow}  \,\,
k[\Lambda]^{36} \,\, {\buildrel \partial \over \longrightarrow}  \,\,
k[\Lambda]^{18} \,\, {\buildrel \partial \over \longrightarrow}  \,\,
k[\Lambda]^9 \,\, {\buildrel \partial \over \longrightarrow}  \,\,
k[\Lambda]^4 \,\, {\buildrel \partial \over \longrightarrow}  \,\,
k $$
One of the $36$ fully attached monomials of degree four is $\, adad $.
It admits three relevant reductions
$ \, adad \,{\buildrel 0 \over \longrightarrow}  \, adad \,$,
$ \, adad \,{\buildrel 1 \over \longrightarrow}  \, cbad \,$  and
$ \, adad \,{\buildrel 3 \over \longrightarrow}  \, bbdb $.
This implies
$$ \partial (E_{adad} ) \,\, = \,\,
  a \cdot E_{dad} \,-\, c \cdot E_{bad} \,-\, b \cdot E_{bdb} .$$

\vskip .3cm
\noindent {\em Example 2} {\sl (The Koszul complex): }
Let $\Lambda = \N^d$. The Gr\"obner basis $\Gc$
consists of the relations $\,y_i y_j \rightarrow y_j y_i \,$
for $\, 1 \leq j \! < \! i \leq n$. A word $w$ is fully attached
if and only if $w = y_{i_1} y_{i_2} \cdots y_{i_r}$
for $\,i_1 \! >\! i_2 \!>\! \cdots \!>\! i_r $. In this case
$\,\partial( E_w ) \,= \,
\sum_{ j=1}^r (-1)^{r-j} y_{i_j} E_{w_j} \,$
where $\,w_j =
y_{i_1} \cdots
y_{i_{j-1}}
y_{i_{j+1}} \cdots
y_{i_r} $. Hence $\, ({\bf F}, \partial) \,$
is the Koszul complex on $n$ indeterminates.
\end{Ex}
\begin{Ex}[The Cartan-complex]\label{aussere} If $A$ is the exterior algebra, then $F_\bullet$ with
\begin{eqnarray*}
F_i&:=&\bigoplus_{\begin{array}{c}1\le j_1<\ldots <j_r\\l_1\ldots l_r\in\N\\i=\sum_{t=1}^rl_t\end{array}}A\:e_{i_1}^{(l_1)}\ldots e_{i_r}^{(l_r)}\\
e_{i_1}^{(l_1)}\ldots e_{i_r}^{(l_r)}&\to&\sum_{t=1}^rx_{i_t}\,e_{i_1}^{(l_1)}\ldots e_{i_t}^{(l_t-1)}\ldots e_{i_r}^{(l_r)}
\end{eqnarray*}
defines a minimal resolution of $k$ as $A$-module, called the Cartan-complex.\\[+3mm]
{\em Proof.} For the exterior algebra $A=k(x_1,\ldots, x_n)/\langle x_ix_j+x_jx_i\rangle$
the resolution $F_\bullet$ is by Corollary \ref{red} minimal. The set of reduction rules is given by
$\mathcal{R}:=\{x_i^2\to 0, x_ix_j\stackrel{-1}{\tol} x_jx_i\mbox{ for }i<j\}$.
Then the fully attached tuple are exactly the words
\[(x_{i_1},\ldots, x_{i_1}, x_{i_2},\ldots,  x_{i_2},\ldots x_{i_r},\ldots, x_{i_r})\mbox{, with $1\le i_1<\ldots <i_r\le n$}.\]
Since $x_ix_j$ is reduced to $-x_jx_i$, if $i\neq j$, and each reduction has factor $(-1)$,
we got for each reduction the coefficient $(-1)(-1)=1$. Since $x_ix_i$ is reduced to $0$,
the differential follows.\qed
\end{Ex}

The following example shows, that even in the case where the Gr\"obner basis is not finite
one can apply our theory:
\begin{Ex}Consider the two-side ideal $\aaf=\langle x^2-xy\rangle$.
By \cite{groebner} there does not exist a finite Gr\"obner basis with respect to degree-lex for $\aaf$. 
One can show, that $\aaf=\langle xy^{n}x-xy^{n+1}\mid n\in\N\rangle$ and that 
$\{ xy^{n}x-xy^{n+1}\mid n\in\N\}$ is an infinite Gr\"obner basis with respect
to degree-lex.

If one applies our matching from Lemma \ref{match2} it is easy to see, that the critical cells are
given by tuples of the form
\[[x|y^{n_1}|x|y^{n_2}|x|\ldots|x|y^{n_l}|x]~\mbox{and}~[x|y^{n_1}|x|y^{n_2}|x|\ldots|x|y^{n_l}],\]
with $n_1,\ldots, n_l\in\N$.\\
By degree-reasons it follows, that the Morse complex is even a minimal resolution.
Therefore we get a minimal resolution $F_\bullet$ of $k$ over
$A=k\langle x_1,\ldots x_n\rangle/\aaf$.

In this case, this proves that $k$ does not admit a linear
resolution and hence $A$ is not Koszul.
\end{Ex}

%------------------------------------------------------------------------------
%
% Application to Hochschild homology
%
%------------------------------------------------------------------------------
\section{application to the acyclic Hochschild complex}\label{hochcase}
Hochschild homology is defined for arbitrary rings. Since we only consider the
case $A=k\langle x_1,\ldots, x_n\rangle/\aaf$, we introduce the Hochschild
homology only for $k$-algebras of this type. For the general definition see
\cite{hoch}. Our proofs still holds, if $A$ is an $R$-algebra, where $R$ is
a commutative Ring and $A$ is free as an $R$-module.

Let $A=k\langle x_1,\ldots, x_n\rangle/\aaf$ be the quotient algebra of a
non-commutative polynomial ring by a two-side ideal $\aaf$. The Hochschild
complex $\HComplex^A_\bullet = (HC_i,\partial_i)$ is defined as follows:
\begin{eqnarray*}
HC_i&:=&A^{\tensor_k i+2},\mbox{ for }i\le -1,\\
\partial(a_0\tensor\ldots\tensor a_{i+1})&:=&
   \sum_{j=0}^{i}(-1)^ja_0\tensor\ldots\tensor a_ja_{j+1}\tensor
                             \ldots\tensor a_{i+1}.
\end{eqnarray*}
Here $A\tensor A^{\op}$ acts on $HC_i$ via
\[(\mu\tensor \gamma)(a_0\tensor\ldots\tensor a_{i+1})=(\mu a_0)\tensor
  a_1\tensor\ldots\tensor a_i\tensor (a_{i+1})\gamma.\]
Let $W$ be a basis of $A$ as $k$-vectorspace, such that $1\in W$. Then we can
fix a basis of $HC_i$ as $(A\tensor_k A^{\op})$-module:
\[\Big\{[w_1|\ldots|w_i]\mbox{, with }w_j\in W.\Big\}\]
Here we identify
$[w_1|\ldots|w_i]$ with $1\tensor w_1\tensor\ldots\tensor w_i\tensor 1$.

In \cite{hoch} it is proved, that $\HComplex^A_\bullet$ is a projective
resolution of $A$ as an $A\tensor A^{\op}$-module,
called the standard Hochschild resolution, or the acyclic Hoch\-schild complex.

We now consider the normalized acyclic Hochschild complex:

The following complex is called normalized standard Hochschild resolution,
which is homotopic to the standard Hochschild resolution
(a proof for the homotopy equivalence via algebraic Morse theory is given in
the appendix):
\begin{enumerate}
  \item $HC_i^{\mathcal{M}}$ has basis
    $\{[w_1|\ldots|w_i]\mid w_j\in W\setminus\{1\}\}.$
  \item The differential is given by
    \begin{eqnarray*}
      \partial_i([w_1|\ldots|w_i])= (w_1\tensor 1)\:[w_2|\ldots|w_i]+
              (-1)^i(1\tensor w_i)\:[w_1|\ldots|w_{i-1}]\\
        +\sum_{j=1}^{i-1}(-1)^j\left(\sum_l a_l\:[w_1|\ldots
              |w_{j-1}|w'_l|w_{j+1}\ldots|w_i]\right),
    \end{eqnarray*}
    if $w_jw_{j+1}=a_0+\sum_l a_l\: w'_l$, with $a_0, a_l\in k$ and
    $w_l'\in W\setminus\{1\}$.
\end{enumerate}
If $M$ is a $A$-bimodule, we regard it as a right $A\tensor A^{\op}$-module
via $m(\mu\tensor \gamma):=\gamma m \mu$.\\
The Hochschild homology $HH(A,M)$ of $A$ with coefficients in $M$ is defined
by the homology of the complex
$M\tensor_{(A\tensor A^{\op})}\HComplex^A_\bullet$.
In this case it follows, that
$HH(A,M)\iso \Tor_\bullet^{A\tensor A^{\op}}(M,A)$ (see \cite{hoch}).

\medskip

%-----------------------------------------------------------------------
%end of general
%-----------------------------------------------------------------------
Now, let $A=k\langle x_1,\ldots,x_n\rangle/\langle f_1,\ldots,f_s\rangle$ be
the non-commutative (resp. commutative) polynomial ring in $n$ indeterminates
divided by a two-side ideal, where $f_1,\ldots, f_s$ is a finite reduced
Gr\"obner basis of $\aaf=\langle f_1,\ldots,f_s\rangle$ with respect to the
degree-lex order. We now give an acyclic matching on the acyclic Hochschild
complex, which is minimal in special cases. Let $\mathcal{G}$ be the set of
standard monomials of degree $\ge 1$ with respect to the degree-lex order.
In this case the normalized acyclic Hochschild complex is given by:
\begin{eqnarray*}
  HC_i&:=&\bigoplus_{w_1,\ldots,w_i\in \mathcal{G}} (A\tensor A^{\op})\;
    [w_1|\ldots|w_i],
\end{eqnarray*}
with differential
\begin{eqnarray*}
  \lefteqn{\partial([w_1|\ldots|w_i]):=(w_1\tensor 1)\,[w_2|\ldots|w_i]}\\
  &&+(-1)^i(1\tensor w_i)[w_1|\ldots|w_{i-1}]\\
  &&+\sum_{j=1}^{i-1}(-1)^{j}\left(\sum_{r}a_r [w_1|\ldots
                   |w_{j-1}|v^j_r|w_{j+2}|\ldots|w_i]\right),
\end{eqnarray*}
if $w_jw_{j+1}$ is reducible to $a_0+\sum_r a_rv^j_r$ (if $w_jw_{j+1}\in\Gc$
we set $v^j_r=w_jw_{j+1}$).\\

We apply the same acyclic matching as in Section \ref{noncommcase}
(resp. Section \ref{commcase}).\\
Since in addition in this case the element $[w_1|\ldots|w_i]$ also maps to 
$(-1)^i(1\tensor w_i)[w_1|\ldots|w_{i-1}]$ we have to modify the
differential:\\
The reduction-rules are the same as in Section \ref{commcase}, except that
the reduction coefficient in Definition \ref{third} is $(c\tensor 1)$ instead
of $c$. In order to define the coefficient we say $e$ can be reduced to $f$
with coefficient $c$ (we write $e\stackrel{c}{\tol}f$), where
$e=(w_1,\ldots, w_l)$ and $f=(v_1,\ldots, v_{l-1})$ are two fully attached
$l$ (resp. $l-1$)-tuples, if there exists a sequence of $l$-tuples
$e=e_0, e_1,\ldots,e_{r-1}$ such that either there exists:
\begin{enumerate}
\item an $l$-tuple $e_r=(u,f)$, such that $e_0$ can reduced to $e_r$ with
  reductions of Type I and III, i.e.
  \[e_0\stackrel{-a_1}{\tol} e_1\stackrel{-a_2}{\tol}e_2\stackrel{-a_3}{\tol}
       \ldots\stackrel{-a_r}{\tol}e_r;\]
  in this case we set $c:=\left((-1)^r\prod_{i=1}^ra_i\right)\,(u\tensor 1)$,
  or
\item an $l$-tuple $e_r=(f,u)$ such that $e_0$ can reduced to $e_r$ with
  reductions of Type I and III, i.e.
  \[e_0\stackrel{-a_1}{\tol} e_1\stackrel{-a_2}{\tol}e_2\stackrel{-a_3}{\tol}
      \ldots\stackrel{-a_r}{\tol}e_r;\]
  in this case we set
  $c:=\left((-1)^{r+k}\prod_{i=1}^ra_i\right)\,(1\tensor u)$, or
\item an $l$-tuple $e_r$, such that $e_0$ can reduced to $e_r$ with reductions
  of Type I and III and $e_r$ can reduced to $f$ with a reduction of Type II, i.e.
  \[e_0\stackrel{-a_1}{\tol} e_1\stackrel{-a_2}{\tol}e_2\stackrel{-a_3}{\tol}
  \ldots\stackrel{-a_r}{\tol}e_r\stackrel{(-1)^jb}{\tol} f;\]
  in this case we set $c:=(-1)^{r+j}\,b\,\prod_{i=1}^ra_i$.
\end{enumerate}
We define the reduction coefficient $[e,f]$ and the complex $F_\bullet$ as in
Section \ref{noncommcase} (resp. Section \ref{commcase}).
With the same proof as in Section \ref{noncommcase} (resp. Section
\ref{commcase}) we obtain the following theorem:
\begin{Thm}
$(F_\bullet, \partial)$ is a free resolution of $A$ as an
$A\tensor A^{\op}$-module.\\
If no reduction of Type II is possible, then $(F_\bullet,\partial)$ is
minimal.\qed
\end{Thm}
Moreover we get similar results to the results from Section \ref{commcase}, \ref{noncommcase}
about minimality of $F_\bullet$ and rationality of
the Poincar\'e-Betti series
\[\poinr{k}{A\tensor A^{\op}}=\sum_{i,\alpha}\dim_k\left((\Tor^{(A\tensor A^{\op})}_i(k,A))_\alpha\right)
     x^\alpha t^i\] from Section \ref{noncommcase}
(resp. Section \ref{commcase}).\\
As in Section \ref{commcase} we can give an explicit description of the
minimal resolution $F_\bullet$ in the following cases:
\begin{enumerate}
\item $A=S/\langle f_1,\ldots, f_s\rangle$, where $S=k[x_1,\ldots, x_n]$ is
  the commutative polynomial ring in $n$ indeterminates and $f_i$ a reduced
  Gr\"obner basis with respect to the degree-lex order, such that the
  initial ideal is a complete intersection (note, in case $s=1$ this
  resolution was first given by BACH (see. \cite{bach})).
\item $A=E$, where $E$ is the exterior algebra.
\end{enumerate}
Let $A=k[x_1,\ldots, x_n]/\langle f_1,\ldots, f_s\rangle$ be the commutative
polynomial-ring in $n$ indeterminates, with
$f_i=x^{\gamma_i}+\sum_{\alpha_i\neq 0}f_{i,\alpha_i}x^{\alpha_i}$, $1\le i\le s$ is a
reduced Gr\"obner basis with respect to the degree-lex order, such that
$x^{\gamma_i}$ is the leading term (since we start with the normalized
Hochschild resolution, the condition $\alpha\neq 0$ is no restriction).

Let $\Gc=\{x^\alpha\mid x^\alpha\not\in \langle x^{\gamma_1},\ldots,
x^{\gamma_s}\rangle\}$ be the set of standard monomials of degree $\ge 1$.
We assume that the initial ideal
$\ini(\aaf)=\langle x^{\gamma_1},\ldots, x^{\gamma_s}\rangle$ is a complete
intersection. With the same arguments as in Theorem \ref{bachkomm} 
it follows that $F_\bullet$ is minimal.
%To give an explicit description of the differential we assume further that
%$x^{\gamma_i}$ is reducible to
%$-\sum_{\alpha\in\mathbb{N}^n}f_{i\alpha}x^{\alpha}$.
% and for
%$j\neq i$, $f_{i\alpha}\neq 0$ and $1\le r\le n$ we have
%$x^{\gamma_j}\,{\not |}\,x_r\frac{x^{\alpha}}{x_{\alpha}}$, where
%$x_\alpha:=\max_{\lex}\big\{x_j~\mbox{divides}~x^\alpha\big\}$.\\
We use the same notations as \cite{bach} and write
\begin{eqnarray*}
  T(x_i)&=&(x_i\tensor 1)-(1\tensor x_i),\\
  \frac{T_i(f)}{T(x_i)}&=&\sum_{\alpha\in\mathbb{N}^n}f_\alpha
       \sum_{j=0}^{\alpha_i-1}(x^{\alpha_1}
  \cdots x^{\alpha_{i-1}}x^j\tensor x^{\alpha_i-1-j}x^{\alpha_{i+1}}\cdots
      x^{\alpha_n}).
\end{eqnarray*}
Under these conditions we get the following result:
%
%theorem
%
\begin{Thm}\label{bachver}Let $A=S/\langle f_1,\ldots, f_s\rangle$, such that the initial
ideal $\ini(\langle f_1,\ldots, f_s\rangle)$ is a complete intersection. 
Then the following complex is a multigraded minimal
resolution of $A$ as an $A\tensor A^{\op}$-module and carries the structure
of a differential graded algebra.
  \[F_i:=\bigoplus_{\begin{array}{c}1\le i_1<\ldots<i_r\le n\\
       1\le j_1<\ldots<j_q\le s\\l_1,\ldots, l_q\in\mathbb{N}\\
       i=r+2\sum_{j=1}^ql_j\end{array}}
   A\tensor A^{\op}\; e_{i_r}\ldots e_{i_1}t_{j_1}^{(l_1)}\ldots
   t_{j_q}^{(l_q)},\]
  \begin{eqnarray*}
    e_{i_r}\ldots e_{i_1}&\mapsto&\sum_{m=1}^r (-1)^{\#\{i_j>i_m\}}T(x_{i_m})
    \,e_{i_r}\ldots \widehat{e_{i_m}}\ldots e_{i_1},\\
    t_{j_1}^{(l_1)}\ldots t_{j_q}^{(l_q)}&\mapsto&
    \sum_{m=1}^q\sum_{p=1}^n
    \frac{T_p(f_{j_m})}{T(x_p)}\;e_pt_{j_1}^{(l_1)}\ldots
    t_{j_m}^{(l_{j_m}-1)}\ldots t_{j_q}^{(l_q)},
  \end{eqnarray*}
  where $t_{i_j}^{(0)}:=1$, $e_ie_j=-e_je_i$ and $e_ie_i=0$. For the
  differential we have:
  \begin{eqnarray*}
    \partial(e_{i_r}\ldots e_{i_1}t_{j_1}^{(l_1)}\ldots t_{j_q}^{(l_q)})&=&
    \partial(e_{i_r}\ldots e_{i_1})t_{j_1}^{(l_1)}\ldots t_{j_q}^{(l_q)}\\
    &&+(-1)^r\,e_{i_r}\ldots e_{i_1}\partial(t_{j_1}^{(l_1)}\ldots
    t_{j_q}^{(l_q)}).
  \end{eqnarray*}
\end{Thm}
Note, that in case $A=S/\langle f\rangle$ this result was first obtained in \cite{bach}
and our complex coincides with the complex given in \cite{bach}.
%
%corollaries
%
\begin{Cor}Under the assumptions of Theorem \ref{bachver} the Hilbert series of the 
Hochschild homology of $A$ with coefficients in $k$
has the form:
\begin{eqnarray*}
  \hilb{HH(A,k)}&=&\sum_{i,\alpha}\dim_k\left(
    (\Tor_i^{A\tensor A^{\op}}(k,A)_\alpha\right)x^\alpha\,t^i\\
  &=&\frac{\D \prod_{i=1}^n(1+x_i\,t)}{\prod_{i=1}^k(1-x^{\gamma_i}t^2)}.
\end{eqnarray*}
\end{Cor}
If $\aaf$ is the zero-ideal, we get with the same arguments the following
special case:
\begin{Cor}
Let $A=k[x_1,\ldots, x_n]$, then the following complex is a minimal resolution
of $A$ as an $A\tensor A^{\op}$-module.
\begin{eqnarray*}
  F_i&:=&\bigoplus_{1\le i_1<\ldots<i_r\le n}A\tensor A^{\op}\;
  e_{i_1}\ldots e_{i_r}\\
  e_{i_1}\ldots e_{i_r}&\mapsto&\sum_{m=1}^r (-1)^{\#\{i_j<i_m\}}T(x_{i_m})\,
  e_{i_1}\ldots \widehat{e_{i_m}}\ldots e_{i_r}
\end{eqnarray*}
In particular we have:
\begin{eqnarray*}
  \hilb{HH(A,k)}&=&\sum_{i,\alpha}\dim_k\left(
    (\Tor_i^{A\tensor A^{\op}})(k,A)_\alpha\right)x^\alpha\,t^i\\
  &=&\prod_{i=1}^n(1+x_i\,t).
\end{eqnarray*}
\end{Cor}
%
%proof of theorem
%
\begin{proof}[Proof of Theorem \ref{bachver}] The description of the basis of $F_i$ 
follows with exactly the same arguments as for the proof of Theorem \ref{bachkomm}.
Since no constant term appears in the differential it suffices to verify that the 
differential has the given form.

First we consider a variable $[x_i]$. Clearly it maps to $(x_i\otimes 1)-(1\otimes x_i)$.

Next we consider a minimal fully attached tuple 
$w_j=\left|x_{\gamma_i}\left|\frac{x^{\gamma_i}}{x_{\gamma_i}}\right.\right|$,
where $x_\gamma:=x_{m(x^\gamma)}$.
Then we have:
\[\left|x_{\gamma_i}\left|\frac{x^{\gamma_i}}{x_{\gamma_i}}\right.\right|\to
    -\sum_{\alpha}f_{i\alpha}|x^{\alpha}|
    \to \sum_{\alpha}f_{i\alpha}\left|x_{\alpha}\left|
    \frac{x^{\alpha}}{x_{\alpha}}\right.\right|.\]
As in the commutative case the multi-index $\alpha$ decreases successively,
but here $\left[x_{\beta}\left|\frac{x^{\alpha}}{x_{\beta}x_{\alpha}}\right.\right]$,
for $x_\beta=x_{m(x^\beta)}$ with $x^\beta:=\frac{x^{\alpha}}{x_{\alpha}}$,
maps in addition to
$\left(1\tensor \frac{x^{\alpha}}{x_{a'}x_{\alpha}}\right)[x_{a'}]$ hence in
this case we get:
\[\left|x_{\gamma_i}\left|\frac{x^{\gamma_i}}{x_{\gamma_i}}\right.\right|\to
\sum_{j=1}^n\frac{T_j(f_i)}{T(x_j)}\;\;e_j.\]

For a fully attached tuple $[w_1,\ldots, w_l]$, we have to calculate the sign of the
permutations. This calculation is similar to the calculation of the sign in the 
commutative case (see proof of Theorem \ref{bachkomm}) and is left to the reader.

With the bijection between the elements $e_{i_r}\ldots e_{i_1}t_{j_1}^{(l_1)}\ldots t_{j_q}^{(l_q)}$ 
and the fully attached tuples we finally get the following differential:
\begin{eqnarray*}
  \lefteqn{\partial(e_{i_r}\ldots e_{i_1}t_{j_1}^{(l_1)}\ldots
      t_{j_q}^{(l_q)})}\\
  &=&\sum_{m=1}^r (-1)^{\#\{i_j>i_m\}}T(x_{i_m})\,
    e_{i_r}\ldots \widehat{e_{i_m}}
    \ldots e_{i_1}t_{j_1}^{(l_1)}\ldots t_{j_q}^{(l_q)}\\
  &&+\sum_{m=1}^q\sum_{p=1\atop p\neq i_1,\ldots,i_r}^n(-1)^r\:
  \frac{T_p(f_{j_m})}{T(x_p)}\;e_{i_r}\cdots e_{i_1}e_pt_{j_1}^{(l_1)}\ldots
   t_{j_m}^{(l_{j_m}-1)}\ldots t_{j_q}^{(l_q)},
\end{eqnarray*}
and the desired result follows.
\end{proof}
%
%
%case exterior algebra:
%
%
We now consider the exterior algebra:
\begin{Thm}\label{bachvernonkomm} Let
$E=k[x_1,\ldots, x_n]/\langle x_i^2, x_ix_j+x_jx_i\rangle$ be the exterior
algebra. The following complex is a minimal resolution of $E$ as
$E\tensor E^{\op}$-module:
  \[F_i:=\bigoplus_{\begin{array}{c}1\le i_1<\ldots<i_r\le n\\
     l_1,\ldots, l_r\in\mathbb{N}^n\end{array}}
  E\tensor E^{\op}\; e_{i_1}^{(l_1)}\ldots e_{i_r}^{(l_r)},\]
  with
  \[e_{i_1}^{(l_1)}\ldots e_{i_r}^{(l_r)}\mapsto
     \sum_{j=1}^r (x_{i_j}\tensor 1)+(1\tensor x_{i_j})\, e_{i_1}^{(l_1)}
          \ldots e_{i_j}^{(l_j-1)}\ldots e_{i_r}^{(l_r)}.\]
  In particular we have :
\[\hilb{HH(E,k)}=\sum_{i,\alpha}\dim_k\left(
   (\Tor_i^{E\tensor E}(k,E))_\alpha\right)x^\alpha\,t^i.\]
Let $S$ be the commutative polynomial-ring in $n$ indeterminates then we have
the following duality:
\begin{eqnarray*}
  \hilb{HH(E,k)}&=&\hilb{S},\\
  \hilb{HH(S,k)}&=&\hilb{E}.
\end{eqnarray*}
\end{Thm}
\begin{proof}
  The proof is the same as in Example \ref{aussere} from Section
  \ref{noncommcase}, with the modified differential.
\end{proof}

%------------------------------------------------------------------------------
%
% Appendix
%
%------------------------------------------------------------------------------

\appendix

\section{The Bar- and the Hochschild-complex}
\label{appendixa}

In this section, we show how to obtain the normalized Bar-resolution
(resp. the normalized acyclic Hochschild complex) via algebraic discrete
Morse theory from the Bar resolution (resp. the acyclic Hochschild complex).
Again we consider the Bar resolution (resp. the acyclic Hochschild complex)
only for $k$-algebras $A$, but the proofs still holds for $R$-algebras $A$,
where $R$ is a commutative Ring and $A$ is projective as an $R$-module.
For the general definition of the Bar resolution (resp. the acyclic Hochschild
complex) see \cite{weib} chap. 8.1. (resp. \cite{hoch}, chap. 2.11.).\\

Let $A$ be a $k$-algebra and let $W$ be a basis of $A$ as a $k$-vectorspace
such that $1 \in W$ and $M$ an $A$-module. The Bar-resolution 
$$\BarRes_\bullet : \cdots \mathop{\rightarrow}^{\partial_{i+1}}
B_i \mathop{\rightarrow}^{\partial_i} B_{i-1}
\mathop{\rightarrow}^{\partial_{i-1}}\cdots \mathop{\rightarrow}^{\partial_1}
B_1 \mathop{\rightarrow}^{\partial_1} B_0 = k$$
of $M$ with respect to $W$ is defined by:
\begin{enumerate}
\item $B_i$ is the free $A\tensor_k M$-module, generated by the tuples
  $[w_1|\cdots|w_i]$, with $w_1,\ldots,w_i \in W$.
\item The differential $\partial_i : B_i \rightarrow B_{i-1}$ is given by
  \begin{eqnarray*}
    \lefteqn{\partial_i ([w_1|\cdots|w_i])=(w_1\otimes 1)~[w_2|\ldots|w_i]}\\
    &&+\D{\sum_{j=1}^{i-1}}(-1)^j \left(\begin{array}{c}
        (a_0\otimes 1)~[w_1|\ldots|w_{j-1}|1|w_{j+2}|\ldots|w_i]\\[+2mm]
        +\sum_l (a_l\tensor 1)~[w_1|\ldots|w_{j-1}|w'_l|w_{j+1}\ldots|w_i]
    \end{array}\right)\\
    &&+(-1)^i~(1\tensor w_i)~[w_1|\ldots|w_{i-1}],
  \end{eqnarray*}
  if $w_jw_{j+1}=a_0+\sum_l a_l\: w'_l$, with $a_0, a_l\in k$ and
  $w_l'\in W\setminus\{1\}$.
\end{enumerate}

\begin{Prop}[Normalized Bar-Resolution]\label{normbar}
  Let $A$ be a $k$-algebra and let $W$ be a basis of $A$ as a $k$-vectorspace
  such that $1 \in W$ and $M$ an $A$-module. Then there is an acyclic matching
  $\Mmc$ on the Bar resolution $\BarRes_\bullet$ with respect
  to $W$ such that the corresponding Morse-complex $\BarRes^{\Mmc}_\bullet$ is
  given by:
  \begin{enumerate}
    \item $B_i^{\Mmc}$ is the free $A\tensor M$-module, generated by the
      tuples $[w_1|\cdots|w_i]$, with $w_1,\ldots,w_i \in W\setminus\{1\}$.
    \item The Mores-differential $\partial^{\Mmc}_i$ is given by
      \begin{eqnarray*}
        \lefteqn{\partial^{\Mmc}([w_1|\ldots|w_i])=
          (w_1\tensor 1)~[w_2|\ldots|w_i]}\\
         &&+\D{\sum_{j=1}^{i-1}}(-1)^j \sum_l (a_l\tensor 1)~
         [w_1|\ldots|w_{j-1}|w'_l|w_{j+1}\ldots|w_i]\\
         &&+(-1)^i~(1\tensor w_i)~[w_1|\ldots|w_{i-1}],
       \end{eqnarray*}
       if $w_jw_{j+1}=a_0+\sum_l a_l\: w'_l$, with $a_0, a_l\in k$ and
       $w_l'\in W\setminus\{1\}$.
  \end{enumerate}
  In particular, $\BarRes^\Mmc_\bullet$ is the normalized bar-resolution.
\end{Prop}
\begin{proof} We define the matching ${\Mmc}$ by
  \[[w_1|\ldots|w_l|w_{l+1}|\ldots|w_i]\to
       [w_1|\ldots|w_lw_{l+1}|\ldots|w_i]\in\Mmc,\]
  if $w_l:=\min(j\mid w_j=1)$,
  $w_{l'}:=\max(j\mid w_r=1\mbox{ for all $l\le r\le j$})$ and $l'<i$ and $l'-l$ is
  odd. The invertibility is given, since in both cases the coefficient in the
  differential is $\pm 1$:
  $$\pm 1\:[w_1|\ldots|w_lw_{l+1}|\ldots|w_i]\in\partial([w_1|\ldots|w_l|w_{l+1}|\ldots|w_i].$$
  It is easy to see, that the other conditions of an acyclic matching are
  satisfied as well. The critical cells are exactly the desired basis
  elements and an element $[w_1|\ldots|w_i]$, for which $w_j=1$ for some $j$
  is never mapped to an element $[w_1|\ldots|w_i]$, with $w_j\neq 1$ for all
  $j$. This implies the formula for the Morse-differential.
\end{proof}

Let $A$ be a $k$-algebra and let $W$ be a basis of $A$ as a $k$-vectorspace
such that $1 \in W$. The acyclic Hochschild-complex
$$\HComplex^A_\bullet : \cdots \mathop{\rightarrow}^{\partial_{i+1}}
   C_i \mathop{\rightarrow}^{\partial_i} C_{i-1}
   \mathop{\rightarrow}^{\partial_{i-1}}
   \cdots \mathop{\rightarrow}^{\partial_1} C_1
   \mathop{\rightarrow}^{\partial_1} C_0 = k$$
with respect to $W$ defined by:
\begin{enumerate}
  \item $C_i$ is the free $(A\tensor_k A^{\op})$-module generated by
    $[w_1|\cdots|w_i]$, with $w_1,\ldots,w_i \in W$.
  \item The differential $\partial_i$ is given by
    \begin{eqnarray*}
      \partial_i([w_1|\ldots|w_i])= (w_1\tensor 1)\:[w_2|\ldots|w_i]+
          (-1)^i(1\tensor w_i)\:[w_1|\ldots|w_{i-1}]\\
      +\sum_{j=1}^{i-1}(-1)^j\left(\begin{array}{c}a_0\:
          [w_1|\ldots|w_{j-1}|1|w_{j+2}|\ldots|w_i]\\[+2mm]
      +\sum_l a_l\:[w_1|\ldots|w_{j-1}|w'_l|w_{j+1}\ldots|w_i]
    \end{array}\right),
    \end{eqnarray*}
    if $w_jw_{j+1}=a_0+\sum_l a_l\: w'_l$, with $a_0, a_l\in k$ and
    $w_l'\in W\setminus\{1\}$.
\end{enumerate}
\begin{Prop}[Normalized Acyclic Hochschild Complex]\label{normhoch}
  Let $A$ be a $k$-Algebra and let $W$ be a basis of $A$ as a
  $k$-vectorspace such that $1 \in W$. Then there is an acyclic matching
  $\Mmc$ on the Hochschild-complex $\HComplex^A_\bullet$ of $A$ such that the
  corresponding Morse-complex $\HComplex^{\Mmc}_\bullet$ is given by:
  \begin{enumerate}
    \item $C_i^{\Mmc}$ is the free $(A\tensor_k A^{\op})$-module generated by
      $[w_1|\cdots|w_i]$, with $w_1,\ldots,w_i \in W\setminus\{1\}$.
    \item The Morse-differential $\partial^{\Mmc}_i$ is given by
      \begin{eqnarray*}
        \partial_i([w_1|\ldots|w_i])= (w_1\tensor 1)\:[w_2|\ldots|w_i]+
            (-1)^i(1\tensor w_i)\:[w_1|\ldots|w_{i-1}]\\
        +\sum_{j=1}^{i-1}(-1)^j\left(\sum_l a_l\:
          [w_1|\ldots|w_{j-1}|w'_l|w_{j+1}\ldots|w_i]\right),
      \end{eqnarray*}
      if $w_jw_{j+1}=a_0+\sum_l a_l\: w'_l$, with $a_0, a_l\in k$ and
      $w_l'\in W\setminus\{1\}$.
  \end{enumerate}
\end{Prop}
\begin{proof} The proof is essentially identical to the proof of Proposition \ref{normbar}.
\end{proof}

%---------------------------------------------------------------
%
%section: Algebraic Discrete Morse Theory
%
%---------------------------------------------------------------
\section{Algebraic Discrete Morse Theory}\label{proof}

In this section we give our prove of the Algebraic Discrete Morse theory
(Theorem \ref{morse}). We write $\Gd(c,c')$ (resp. $\Gu(c,c')$) for the sum
of the weights of all those paths from $c$ to $c'$, for which the first
step $c\to c_1$ satisfies $c\in X_i^\Mmc$ and $c_1\in X^\Mmc_{i-1}$
(resp. $c'\in X^\Mmc_{i+1}$) In most cases it will be clear from the context,
e.g. if $c$ is critical, whether the first step increases or decreases
dimension. Still for the sake of readability we will always equip $\Gamma$
with the respective arrow.

>From now on we assume always, that $\Mmc$ satisfy the three conditions.

We first prove, that the Morse-differential satisfies
$\partial^{\Mmc}_i\circ\partial^{\Mmc}_{i+1}=0$.

%---------------------------------------------------------------
%
%Lemma \partial^2=0:
%
%---------------------------------------------------------------
\begin{Lem}\label{rauf=runter}
Let $\Mmc\subset E$ be an acyclic matching on $G(C_\bullet)=(V,E)$. Then
\begin{itemize}
\item[(P1)] $\partial^{\Mmc}$ is a differential, i
  (i.e. $\partial^{\Mmc}\circ\partial^{\Mmc}=0$).
\item[(P2)] For $(\alpha,\beta,[\alpha:\beta])\in\Mmc$, with
  $\alpha\in X_{i+1}$, $\beta\in X_i$ we have for all $c\in X^\Mmc_{i-1}$
  \[\Gd(\beta,c)=\sum_{c'\in X_i^\Mmc}\Gu(\beta,c')\Gd(c',c))\]
\end{itemize}
\end{Lem}
\disbreak
%---------------------------------------------------------------
%
%Proof of Lemma \partial^2=0:
%
%---------------------------------------------------------------
\begin{proof}
  The proof is by induction over the cardinality of $\Mmc$. In order to
  prove the induction we assume, that both properties are satisfied for
  smaller matchings.

  Let $\Mmc=\{(\alpha,\beta,[\alpha:\beta])\}$ be a matching of
  cardinality $1$.\\

  \noindent {\sf Property (P2):}
  \begin{eqnarray*}
    \lefteqn{0 = \partial^2(\alpha) =
      \sum_{c'\in X_i^\Mmc}[\alpha:c']\partial(c')+[\alpha:\beta]\partial(\beta)}\\
    &=& \sum_{c\in X_{i-1}}\left(\sum_{c'\in X_i^\Mmc}[\alpha:c'][c':c]\right)c
    +  \sum_{c\in X_{i-1}}[\alpha:\beta][\beta:c]c\\
    &=& -[\alpha:\beta]\sum_{c\in X_{i-1}}
    \left(\sum_{c'\in X_i^\Mmc}\left(-\frac{1}{[\alpha:\beta]}\right)[\alpha:c'][c':c]\right)c\\
    && +  \sum_{c\in X_{i-1}}[\alpha:\beta][\beta:c]c\\
    & = & -[\alpha:\beta]\sum_{c\in X_{i-1}}\sum_{c'\in X_i^\Mmc}\Gu(\beta,c')\Gd(c',c)c
    +  \sum_{c\in X_{i-1}}[\alpha:\beta][\beta:c]c\\
    & = & [\alpha:\beta]\sum_{c\in X_{i-1}}\left([\beta:c]- \sum_{c'\in X_i^\Mmc}\Gu(\beta,c')\Gd(c',c)\right)c\\
    & = & [\alpha:\beta]\sum_{c\in X_{i-1}}\left(\Gd(\beta,c)-
      \sum_{c'\in X_i^\Mmc}\Gu(\beta,c')\Gd(c',c)\right)c
  \end{eqnarray*}
  Since $[\alpha:\beta]\in Z(R)\cap R^*$ is not a zero-divisor, and the critical cells are
  linearly independent we got the desired result:
  \[\Gd(\beta,c)-\sum_{c'\in X_i^\Mmc}\Gu(\beta,c')\Gd(c',c)=0\]
%---------------------------------------------------------------
%
%end of induction for property P2:
%
%---------------------------------------------------------------
  \noindent {\sf Property (P1):} Let $c\in X_{i+1}^\Mmc$ be a critical cell.
  We have to distinguish three cases. Note, that the validity of property
  (P2) has been established above.

  \noindent {\em Case 1:} $(\partial^\Mmc)^2(c)=\partial^2(c)$. Since
  $\partial$ is a differential we are done.

  \noindent {\em Case 2:} There exists elements $\beta\in X_i$ and
  $c\neq\alpha\in X_{i+1}$ with $[c:\beta]\neq 0$ and
  $\{(\alpha,\beta,[\alpha:\beta])\}=\Mmc$. Then we have:
  \begin{eqnarray*}
    \lefteqn{(\partial^\Mmc)^2(c)= \D{\sum_{\beta\neq c'\le c}[c:c']\partial^\Mmc(c')}
    + [c:\beta](-\frac{1}{[\alpha:\beta]})\D{\sum_{c'\in X^\Mmc_i\atop c'\neq \beta}}[\alpha:c']\partial^\Mmc(c')}\\
    &=& \D{\sum_{\beta\neq c'\le c}\sum_{c''\le c'}}[c:c'][c':c'']c'' \\
    && + [c:\beta](-\frac{1}{[\alpha:\beta]})\D{\sum_{c'\in X^\Mmc_{i}\atop c'\neq \beta}
      \sum_{c''\le c'}}[\alpha:c'][c':c'']c''\\
    &=& \D{\sum_{c''\in X_{i-1}^\Mmc}}\left( \D{\sum_{\beta\neq c'\le c}}[c:c'][c':c'']\right.\\
    && \left.+ [c:\beta](-\frac{1}{[\alpha:\beta]})\D{\sum_{c'\in X^\Mmc_{i}\atop c'\neq \beta}}[\alpha:c'][c':c'']\right)c''\\
    &=& \D{\sum_{c''\in X_{i-1}^\Mmc}}\left(
      \D{\sum_{\beta\neq c'\le c}}[c:c'][c':c'']+
      [c:\beta]\D{\sum_{c'\in X^{\Mmc'}_{i}\atop c'\neq \beta}}\Gu(\beta,c')\Gd(c',c'')\right)c''\\
    &\stackrel{(P2)}{=}& \D{\sum_{c''\in X_{i-1}^\Mmc}}\left( \D{\sum_{\beta\neq c'\le c}}[c:c'][c':c'']+ [c:\beta]\Gd(\beta,c'')\right)c''\\
    &=& \sum_{c''\in X_{i-1}}\left(\D{\sum_{c'\le c}}[c:c'][c':c'']\right)c'' = \partial^2(c)=0.
  \end{eqnarray*}

\noindent {\em Case 3:} There exists element $\beta\in X_i$ and
$\alpha\in X_{i-1}$ with $[c:\beta]\neq 0$ and
$\{(\beta,\alpha,[\beta:\alpha])\}=\Mmc$\\
Since $\partial^2(c)=0$ we have
\begin{eqnarray*}
0 & =&  \sum_{c'\le c}[c:c'][c':\alpha]\\
& = & [c:\beta][\beta:\alpha]+\sum_{{c'\le c\atop c'\neq\beta}}[c:c'][c':\alpha]\\
& = & [\beta:\alpha]\left([c:\beta] +
        \sum_{{c'\le c\atop c'\neq\beta}}\frac{1}{[\beta:\alpha]}[c:c'][c':\alpha]\right)\\
\end{eqnarray*}
Since $[\beta:\alpha]\in Z(R)\cap R^*$ is not a zero-divisor it follows
\begin{eqnarray}
[c:\beta] = \sum_{{c'\le c\atop c'\neq\beta}}\left(-\frac{1}{[\beta:\alpha]}\right)[c:c'][c':\alpha].\label{case 3}
\end{eqnarray}
This observation allows us to deduce the desired result:
\begin{eqnarray*}
  \lefteqn{(\partial^\Mmc)^2(c) = \D{\sum_{{c'\le c\atop c'\neq\beta}}}[c:c']\partial^\Mmc(c')}\\
  &=& \D{\sum_{{c'\le c\atop c'\neq\beta}} \sum_{{c''\le c'\atop c''\neq\alpha}}}[c:c'][c':c'']c''\\
  && + \D{\sum_{c''\le\beta\atop c''\neq\alpha}\underbrace{\left(\sum_{{c'\le c\atop c'\neq\beta}}[c:c'][c':\alpha]
         \left(-\frac{1}{[\beta:\alpha]}\right)\right)}_{=[c:\beta]\mbox{ by (\ref{case 3})}}[\beta:c'']c''}\\
  &=& \D{\sum_{{c'\le c\atop c'\neq\beta}} \sum_{{c''\le c'\atop c''\neq\alpha}}}[c:c'][c':c'']c''
  + \D{\sum_{c''\le\beta\atop c''\neq\alpha}[c:\beta][\beta:c'']c''} \\
  &=& \D{\sum_{c'\le c}\sum_{{c''\le c'\atop c''\neq\alpha}}}[c:c'][c':c'']c''=0\mbox{, since }\partial^2=0.
\end{eqnarray*}

%---------------------------------------------------------------
%
%end of induction for property P1:
%
%---------------------------------------------------------------
We now assume properties (P1) and (P2) for matchings of cardinality $\le n$.
Let $\Mmc$ be an acyclic matching of cardinality $n+1$. and
$\Mmc':=\Mmc\setminus \{(\alpha,\beta,[\alpha:\beta])\}$, with
$\alpha\in X^{\Mmc'}_{i+1}$ and $\beta\in X^{\Mmc'}_i$.
Then $\alpha$, $\beta$ are critical with respect to $\Mmc'$, and by induction
$\Mmc'$ satisfies (P1) and (P2).

\noindent {\sf Property (P2):}
\begin{eqnarray*}
0 = (\partial^{\Mmc'})^2(\alpha)
 = &
  \D{\sum_{c'\in X^{\Mmc'}_i}
  \sum_{c\le \alpha}}[\alpha:c]\Gu(c,c')\partial^{\Mmc'}(c')\\
 = & [\alpha:\beta]\partial^{\Mmc'}(\beta)
     + \D{\sum_{c'\in X^{\Mmc'}_i\atop c'\neq \beta}
 \sum_{c\le \alpha}}[\alpha:c]\Gu(c,c')\partial^{\Mmc'}(c')\\
 = & [\alpha:\beta]\partial^{\Mmc'}(\beta)
   +\D \sum_{c'\in X^{\Mmc'}_i\atop c'\neq \beta}\Gd(\alpha,c')\partial^{\Mmc'}(c')\\
 = & [\alpha:\beta]\left(\D\sum_{c'\le \beta}[\beta:c']
   \D{\sum_{c \in X^{\Mmc'}_{i-1}}}\Gu(c',c)c\right) \\
   & + \D{\sum_{c'\in X^{\Mmc'}_i\atop c'\neq\beta}}\Gd(\alpha,c')
       \D{\sum_{c\in X^{\Mmc'}_{i-1}}}\Gd(c',c)c\\
 = & [\alpha:\beta]\left(\D{\sum_{c\in X^{\Mmc'}_{i-1}}}\Gd(\beta,c)c\right) \\
   & - [\alpha:\beta]\D{\sum_{c\in X^{\Mmc'}_{i-1}}}\left(-\frac{1}{[\alpha:\beta]}\right)
  \D{\sum_{c'\in X^{\Mmc'}_i\atop c'\neq\beta}}\Gd(\alpha,c')\Gd(c',c)c\\
 = & [\alpha:\beta]\D{\sum_{c\in X^{\Mmc}_{i-1}}(\Gd(\beta,c) }
   - \D{\sum_{c'\in X^{\Mmc}_i}}\Gu(\beta,c')\Gd(c',c))c
\end{eqnarray*}

Since the critical cells are linearly independent and $[\alpha:\beta]$ is a
unit, we got the desired result:
\[\Gd(\beta,c)- \sum_{c'\in X_{\Mmc}^{(i)}}\Gu(\beta,c')\Gd(c',c)= 0\]
%---------------------------------------------------------------
%
%end of induction for property P2:
%
%---------------------------------------------------------------
\noindent {\sf Property (P1):} Let $c\in X_{i+1}^\Mmc$ be a critical cell.
In order to prove the first statement, we have, as in case $1$, to
distinguish three cases:

\noindent {\it Case 1:} $(\partial^\Mmc)^2(c)=(\partial^{\Mmc'})^2(c)$. Since
by induction $(\partial^{\Mmc'})^2=0$ we are done.

\noindent {\it Case 2:} There exists elements $c\neq\alpha\in X_{i+1}^{\Mmc'}$
and $\beta\in X_i^{\Mmc'}$ with $[c:\beta]\neq 0$ and
$(\alpha,\beta,[\alpha:\beta])\in\Mmc$. Then we have:
\[\partial^\Mmc(c)=
  \sum_{c'\in X^\Mmc_i}[c:\beta]\Gu(\beta,c')c'
  + \sum_{c'\in X^\Mmc_i}\Gd(c,c')c',\]
where the last sum is over all paths, which don't go through $\beta$.
It follows
\begin{eqnarray*}
\lefteqn{ (\partial^\Mmc(c))^2  =  \sum_{c'\in X^\Mmc_i}[c:\beta]\Gu(\beta,c')\partial^\Mmc(c')
  + \sum_{c'\in X^\Mmc_i}\Gd(c,c')\partial^\Mmc(c')}\\
& = & \sum_{c''\in X^\Mmc_{i-1}}[c:\beta]\left(\sum_{c'\in X^{\Mmc}_i}\Gu(\beta,c')\Gd(c',c'')\right)c''
   + \sum_{c''\in X^\Mmc_{i-1}}\Gd(c,c'')c''\\
& = & \sum_{c''\in X^\Mmc_{i-1}}[c:\beta]\Gamma_{\downarrow}(\beta,c'')c''
   + \sum_{c''\in X^\Mmc_{i-1}}\Gd(c,c'')c''\\
& = & \sum_{c''\in X_{i-1}^{\Mmc'}}\Gd(c,c'')c''
= 0\mbox{, since by induction }\partial^{\Mmc'}\circ\partial^{\Mmc'} =0.
\end{eqnarray*}

\noindent {\it Case 3:} There exists elements $\beta\in X_i^{\Mmc'}$ and
$\alpha\in X_{i-1}^{\Mmc'}$ with $[c:\beta]\neq 0$ and
$(\beta,\alpha,[\beta:\alpha])\in\Mmc$. Then we have
\begin{eqnarray}
\nonumber\lefteqn{(\partial^\Mmc)^2(c)=\sum_{c'\in X^\Mmc_i\atop c'\neq\beta}\Gd(c,c')\partial^\Mmc(c')}\\
\nonumber& = & \sum_{c''\neq \alpha}\left(
  \sum_{c'\in X^\Mmc_i\atop c'\neq\beta}\Gd(c,c')[c':\alpha]\left(-\frac{1}{[\beta:\alpha]}\right)\right)\Gd(\beta,c'')c''\\
    && + \sum_{c''\neq\alpha}\left(\sum_{c'\in X^\Mmc_i\atop c'\neq\beta}
   \Gd(c,c')\Gd(c',c'')\right)c''\\\label{only}
\nonumber& \stackrel{(*)}{=} & \sum_{c''\neq \alpha}\Gd(c,\beta)\Gd(\beta,c'')c''
 + \sum_{c''\neq\alpha}\left(\sum_{c'\in X^\Mmc_{i-1}\atop c'\neq\beta}
   \Gd(c,c')\Gd(c',c'')\right)c''\\
\nonumber& = & \sum_{c''\neq\alpha}\left(\sum_{c'\in X_i^{\Mmc'}}
  \Gd(c,c')\Gd(c',c'')\right)c''\\
\nonumber & = & 0\mbox{ ,since }(\partial^{\Mmc'})^2=0.
\end{eqnarray}
(In (\ref{only}) $\Gd(c',c'')$ counts only paths, which don't go through
$\alpha$). In $(*)$ we use the fact
$\D \Gd(c,\beta)=\sum_{c'\in X_i^\Mmc\atop c'\neq\beta}\Gd(c,c')
\left(-\frac{1}{[\beta:\alpha]}\right)[c':\alpha]$, which holds with the same
argument as in (\ref{case 3}).
\end{proof}
%---------------------------------------------------------------
%
%end of proof for lemma \partial^2=0
%
%---------------------------------------------------------------

In the following, we show that the Morse complex is homotopy equivalent to
the original complex. Thereby, it will be possible to minimize a complex of
free $R$-modules by means of Algebraic Discrete Morse theory.

Let $(C(X,R),\partial)$ be a complex of free $R$-modules, $\Mmc\subset E$ a
matching on the associated graph $G(C(X,R))=(V,E)$ and
$(C(X^\Mmc,R),\partial^\Mmc)$ the Morse complex. We consider the following
maps:

\begin{eqnarray}
f:C(X,R)&\to&C(X^\Mmc,R)\\
\nonumber c\in X_i&\mapsto&f(c):=\sum_{c'\in X_i^\Mmc}\Gamma(c,c')c'\\[+2mm]
g:C(X^\Mmc,R)&\to&C(X,R)\\
\nonumber c\in X_i^\Mmc&\mapsto&g_i(c):=\sum_{c'\in X_i}\Gamma(c,c')c'\\[+2mm]
\chi:C(X,R)&\to&C(X,R)\\
\nonumber c\in X_i&\mapsto&\chi_i(c):=\sum_{c'\in X^{(i+1)}}\Gamma(c,c')c'
\end{eqnarray}

Then:
\begin{Lem}\label{1.lemma}
The maps $f$ and $g$ are homomorphisms of complexes of free $R$-modules.
In particular,
\begin{itemize}
\item[(C1)] $\partial^\Mmc\circ f = f\circ\partial$.
\item[(C2)] $\partial\circ g = g\circ\partial^\Mmc$.
\end{itemize}
\end{Lem}

\begin{Lem}\label{2.lemma}
The maps $g$ and $f$ define a chain homotopy. In particular,
\begin{enumerate}
\item[(H1)] $g_i\circ f_i-\id =\partial\circ\chi_{i+1}+\chi_i\circ\partial$,
  i.e. it is null-homotopic,
\item[(H2)] $f_i\circ g_i-\id = 0$, in particular $f \circ g$ is
  null-homotopic.
\end{enumerate}
\end{Lem}
%
%ende lemma2 f,g sind kettenhomotopie
%
\begin{Cor}[Thm. \ref{morse}]
$C(X^\Mmc,R)$ is a complex of free $R$-modules and
\[H_i(C(X,R),R) = H_i(C(X^\Mmc,R),R).\]
\end{Cor}
\begin{proof}
This is an immediate consequence of Lemma \ref{2.lemma}.
\end{proof}

%---------------------------------------------------------------
%
%proof for 1.lemma
%
%---------------------------------------------------------------
\begin{proof}[Proof of Lemma \ref{1.lemma}:]
\noindent {\sf Property (C1):} Let $c\in X_i$. Then:
\[\partial^\Mmc\circ f(c) = \partial^\Mmc\left(
  \sum_{c'\in X_i^\Mmc}\Gu(c,c')c'\right)
 =  \sum_{c''\in X^\Mmc_{i-1}}
  \sum_{c'\in X_i^\Mmc}\Gu(c,c')\Gd(c',c'')c''.
\]
and
\begin{eqnarray*}
f\circ\partial(c) = & f\left(\sum_{c'\le c}[c:c']c'\right) \\
                  = & \D{\sum_{c''\in X^\Mmc_{i-1}}\sum_{c'\le c}}[c:c']\Gu(c',c'')c'' \\
                  = & \D{\sum_{c''\in X^\Mmc_{i-1}}}\Gd(c,c'').\\
\end{eqnarray*}

Using Lemma \ref{rauf=runter} (P2) the assertion now follows.

\medskip

\noindent {\sf Property (C2):} Let $c\in X_i^\Mmc$.
\begin{eqnarray*}
\partial\circ g(c) = & \D{\sum_{c''\le c}[c:c'']c'' +
       \sum_{c'\in X_i}\Gd(c,c')\sum_{c''\le c'}}[c':c'']c''\\
                   = & \underbrace{\D{\sum_{c''\in X^\Mmc_{i-1}}}\Gd(c,c'')c''}_{(A)} \\
                     & + \underbrace{\D{\sum_{c'\in X_i}\Gd(c,c')
       \sum_{c''\le c'\atop (c'',\beta,[c'':\beta])\in\Mmc}}[c':c'']c''}_{(B)}\\
                     & + \underbrace{\D{\sum_{c'\in X_i}\Gd(c,c')
       \sum_{c''\le c'\atop (\beta,c'',[\beta:c''])\in\Mmc}}[c':c'']c''}_{(C)}
\end{eqnarray*}
We have $(C)=0$: Fix $c''\in X_{i-1}$ and $\beta\in X_i$ such that the edge
$(\beta,c'',[\beta:c''])\in\Mmc$. Then:
\begin{eqnarray*}
\D{\sum_{c'\in X_i}}\Gd(c,c')[c':c'']
 = & \D{\sum_{c'\neq\beta}\Gd(c,c')}[c':c'']c'' +\Gd(c,\beta)[\beta:c'']c''\\
 = & \D{\sum_{c'\neq\beta}\Gd(c,c')}[c':c'']c''\\
   & + \left(\Gd(c,c'')\left(-\frac{1}{[\beta:c'']}\right)\right)[\beta:c'']c''\\
 = &\Gd(c,c'')c''-\Gd(c,c'')c''=0.
\end{eqnarray*}

On the other hand:
\begin{eqnarray*}
\lefteqn{g\circ\partial^\Mmc(c)
  = g\left(\sum_{c'\in X^\Mmc_{i-1}}\Gd(c:c')c'\right)}\\
&=& \underbrace{\sum_{c'\in X^\Mmc_{i-1}}\Gd(c,c')c'}_{(A)}
    + \underbrace{\sum_{c'\in X^\Mmc_{i-1}}\sum_{c''\in X_{i-1}\atop (c'',\beta,[c'':\beta])\in\Mmc}\Gd(c,c')\Gd(c',c'')c''}_{(D)}\\
\end{eqnarray*}
We will verify $(B)=(D)$: Consider the matching
$\Mmc'\setminus \{(c'',\beta,[c'',\beta])\}$. Since $c''$ and $\beta$ are
critical cells in $\Mmc'$ it follows by Lemma \ref{rauf=runter} (P1)
(i.e. $(\partial^{\Mmc'})^2=0$) that
\[0=\sum_{c'\in X_{i-1}^{\Mmc'}\atop c'\neq c''}\Gd(c,c')\Gd(c',\beta) +
\Gd(c,c'')[c'':\beta].\]
Multiplying by $\left(-\frac{1}{[c''.\beta]}\right)$ yields:
\[\sum_{c'\in X_{i-1}^{\Mmc}}\Gd(c,c')\Gd(c',c'')
  = \sum_{c'\in X_i}\Gd(c,c')[c':c''].\]
Thus $(B)=(D)$.
\end{proof}
%---------------------------------------------------------------
%
%proof for 2.lemma
%
%---------------------------------------------------------------
\begin{proof}[Proof of Lemma \ref{2.lemma}:]

\noindent {\sf Property (H2):} Let $c\in X_i^\Mmc$. The map $g$ sends $c$ to
a sum over all $c'\in X_i$ that can be reached from $c$. Since $c$ is
critical, $c'$ can be reached from $c$ if either $c = c'$ or there is a
$c''\in X_{i-1}$, such that $(c',c'',[c':c''])\in\Mmc$. Moreover,
\[f(c)=0\mbox{, is there is a $c'\in X_{i-1}$ such that }(c,c',[c:c'])\in\Mmc.\]
Since $f$ and $g$ are $R$-linear it follows that $f_i\circ g_i(c) = f_i(c)$.
>From $f_{X_\Mmc}=\id$ we infer the assertion.

\noindent {\sf Property (H1):}
We distinguish
\noindent {\it Case 1:} Assume $c$ is critical. Then
\[g_i\circ f_i-\id(c) = g_i(c)-c=\sum_{c'\in X_i\atop (c',\beta,[c':\beta])\in\Mmc}\Gd(c,c')c'\]
Moreover, $\chi_i(c)=0$, in particular, $\partial\circ\chi_i(c)=0$.
\begin{eqnarray*}
\chi(\partial(c)) = & \chi\left(\sum_{c'\le c}[c:c']c'\right) \\
                  = & \D{\sum_{c'\le c}[c:c']\sum_{c''\in X_i}}\Gu(c',c'')c''\\
                  = & \D{\sum_{c''\in X_i\atop (c'',\beta,[c':\beta])\in\Mmc}}\Gd(c,c'')c'' = (g_i\circ f_i-\id)(c).
\end{eqnarray*}

\noindent {\it Case 2:} There is an $\alpha\in X_{i-1}$ such that
$(c,\alpha,[c:\alpha])\in\Mmc$. Then $\chi(c)=0$ and $(g_i\circ f_i-\id)(c) = -\id(c) = -c$. Moreover,
\begin{eqnarray*}
\chi(\partial(c)) = & \chi\left(\D{\sum_{c'\le c}}[c:c']c'\right) \\
                  = & \D{\sum_{c'\le c}[c:c']\sum_{c''\in X_i}}\Gu(c',c'')c''\\
                  = & [c:\alpha]\left(-\frac{1}{[c:\alpha]}\right)c \\
                    & + \D{\sum_{c'\le c\atop c'\neq\alpha}[c:c']\sum_{c''\in X_i\atop c''\neq\alpha}}\Gu(c',c'')c''
 + [c:\alpha]\sum_{c''\in X_i}\Gu(\alpha,c'')c''
\end{eqnarray*}
Since
\[\Gu(\alpha,c'')=\left(-\frac{1}{[c:\alpha]}\right)\sum_{c'\le c\atop c'\neq\alpha}[c:c']\Gu(c',c''),\]
the assertion follows.

\noindent {\it Case 3:} There is an $\alpha\in X_{i+1}$ such that $(\alpha,c,[\alpha:c])\in\Mmc$. Then:

\begin{eqnarray*}
(g_i\circ f_i-\id)(c) = & -c  + \underbrace{\D{\sum_{c'\in X_i^\Mmc}}\Gu(c,c')c'}_{(A)} \\
                        & + \underbrace{\D{\sum_{c''\in X_i\atop (c'',\beta,[c'':\beta])\in\Mmc}
                            \sum_{c'\in X_i^\Mmc}}\Gu(c,c')\Gd(c',c'')c''}_{(B)}. \\
\end{eqnarray*}

On the other hand:
\begin{eqnarray*}
\partial\chi(c)= & \partial\left(\D{\sum_{c'\in X_{i+1}}}\Gu(c,c')c'\right) \\
               = & \D{\sum_{c'\in X_{i+1}}\Gu(c,c')\sum_{c''\le c'}}[c':c'']c''\\
               = & \D{\sum_{c''\le \alpha}}\left(-\frac{1}{[\alpha:c]}\right)[\alpha:c'']c'' \\
                 & + \D{\sum_{c'\neq\alpha}\Gu(c,c')\sum_{c''\le c'}}[c':c'']c''\\
               = & -c + (A) + \underbrace{\D{\sum_{c''\in X_i\atop (c'',\beta,[c'':\beta])\in\Mmc}
  \sum_{c'\neq\alpha}\Gu(c,c')\sum_{c''\le c'}}[c':c'']c''}_{(C)}\\
                 & + \underbrace{\D{\sum_{c''\in X_i\atop (\beta,c'',[\beta:c''])\in\Mmc}
  \sum_{c'\neq\alpha}\Gu(c,c')\sum_{c''\le c'}}[c':c'']c''}_{(D)}\\
\end{eqnarray*}
and
\begin{eqnarray*}
\chi\partial(c)=\chi\left(\sum_{c'\le c}[c:c']c'\right) & \\
 = & \D{\sum_{c''\in X_i\atop (c'',\beta,[c'':\beta])\in\Mmc}\sum_{c'\le c}}[c:c']\Gu(c',c'')c''\\
 =&\underbrace{\sum_{c''\in X_i\atop (c'',\beta,[c'':\beta])\in\Mmc}\Gd(c,c'')c''}_{(E)}
\end{eqnarray*}
We show:
\begin{itemize}
\item[(a)] $(D)=0$.
\item[(b)] $(E)+(C)=(B)$.
\end{itemize}

\noindent {\sf Assertion (a);}
Fix $c''\in X_i$ and $\beta\in X_{i+1}$ such that
$(\beta,c'',[\beta:c''])\in\Mmc$. Then:
\begin{eqnarray*}
\sum_{c'\in X_{i+1}}\Gu(c,c')[c':c''] & \\
 =&  \D{\sum_{c'\neq\beta}}\Gu(c,c')[c':c'']c'' +\Gu(c,\beta)[\beta:c'']\\
 =&  \D{\sum_{c'\neq\beta}}\Gu(c,c')[c':c'']c''\\
  & + \left(\Gu(c,c'')\left(-\frac{1}{\beta:c'']}\right)\right)[\beta:c'']c''\\
 =&\Gamma(_\uparrow(c,c'')c''-\Gu(c,c'')c''=0.
\end{eqnarray*}

\noindent {\sf Assertion (b);} Let $c''\in X_i$ and $\beta\in X_{i-1}$ such
that $(c'',\beta,[c'':\beta])\in\Mmc$. Consider the matching
$\Mmc'=\Mmc\setminus\{(c'',\beta,[c'':\beta])\}$. Then by
Lemma \ref{rauf=runter} (P2)
\[\sum_{c'\in X_i^{\Mmc'}}\Gu(c,c')\Gd(c',\beta) =
  \Gd(c,\beta)\]
Since $c''$ is critical with respect to $\Mmc'$ it follows that
\[ \sum_{c'\in X_i^{\Mmc'}\atop c'\neq c''}\Gu(c,c')\Gd(c',\beta) +
        \Gu(c,c'')[c'':\beta] =\Gd(c,\beta).\]
Multiplying the equation with $\left(-\frac{1}{[c'':\beta]}\right)$,
yields
\[\sum_{c'\in X_i^\Mmc}\Gu(c,c')\Gd(c',c'')
  = \sum_{c'\in X_{i+1}}\Gu(c,c')[c':c''] + \Gd(c,c''),
 \]
where paths are taken with respect to the matching $\Mmc$.
Hence $(B)=(C)+(E)$.
\end{proof}

%------------------------------------------------------------------------------
%
% References
%
%------------------------------------------------------------------------------

%------------------------------------------------------------------------------
%
% End of document
%
%------------------------------------------------------------------------------
\end{document}